\documentclass[a4paper,10pt]{article}
\usepackage[T2A]{fontenc} %% Заказ шрифтов Лапко-Ходулева
\usepackage{personal}

\usepackage[cp866]{inputenc} %% Текст в альтернативной кодировке

\usepackage[english,russian]{babel} %% Русские автопереносы и стандарты
\usepackage{amssymb,amsmath,color} %% Подключение AmsLaTeX не помешает

\newtheorem{thm}{\rm Т е о р е м а}[sect]
\newtheorem{lem}{\rm Л е м м а}[sect]
\newtheorem{rem}{\rm З а м е ч а н и е}[sect]
\newtheorem{exa}{\rm П р и м е р}[sect]

\newtheorem{cor}{\rm С л е д с т в и е}[sect]

\begin{document}

\setcounter{secnumdepth}{3} \setcounter{tocdepth}{2}

\renewcommand{\contentsname}{\vskip-1.7cm}\contentsname

\begin{center} {\large \bf ИНДЕКС ПУАНКАРЕ И ПЕРИОДИЧЕСКИЕ РЕШЕНИЯ
ВОЗМУЩЕННЫХ АВТОНОМНЫХ СИСТЕМ\footnote{Основные результаты статьи
доложены на семинаре Отдела дифференциальных уравнений
Математического института им.~В.~А.~Стеклова РАН 22 марта
2006~г.}}
\\ \vskip0.3cm
\it О. Ю. Макаренков  (Воронеж)
\end{center}

\vskip0.0cm

\centerline{СОДЕРЖАНИЕ}

\tableofcontents

\subsection{Введение}\nonumber
\subsubsection{}\vskip-1cm
\addtocounter{sect}{1}

Рассмотрим автономную систему обыкновенных дифференциальных
уравнений
\begin{equation}\label{np_mi}
  \dot x=f(x),\quad x\in\mathbb{R}^2,
\end{equation}
где $f$ -- непрерывно дифференцируемая функция такая, что решение
системы (\ref{np_mi}) с любым начальным условием продолжимо на
интервал $(-\infty,\infty).$  Систему (\ref{np_mi}) будем называть
порождающей. Предположим, что система (\ref{np_mi}) допускает
$T$-периодический цикл $\widetilde{x}.$ В силу автономности, при
любой сдвижке $\theta\in[0,T]$ функция
$$
  {\widetilde{x}^\theta(t)=\widetilde{x}(\theta+t)}
$$
вновь является $T$-периодическим решением системы (\ref{np_mi}).
Предположим, что функция
$g:\mathbb{R}\times\mathbb{R}^2\to\mathbb{R}^2$ $T$-периодична по
первой переменной и рассмотрим возмущенную систему
\begin{equation}\label{ps_mi}
  \dot x=f(x)+\varepsilon g(t,x),\quad x\in\mathbb{R}^2.
\end{equation}
Статья
 посвящена классической задаче, восходящей к
А.~Пуанкаре, о существовании у возмущенной системы (\ref{ps_mi})
$T$-периодических решений $\widetilde{x}_{\varepsilon},$
сходящихся при $\varepsilon\to 0$ к порождающему решению
$\widetilde{x}^{\theta}.$ Пуанкаре принадлежит утверждение о том,
что параметр $\theta=\theta_0$ порождающего решения, для которого
такое существование имеет место, необходимо является нулем
некоторой бифуркационной функции $M.$ Функцию $M$ Пуанкаре
записывал в неявном виде.

Позднее И.~Г.~Малкин \cite{mal} и В.~К.~Мельников \cite{mel}
указали явный вид функции $M$ и при помощи теоремы о неявной
функции в различных ситуациях доказали, что достаточным условием
существования у системы (\ref{ps_mi}) $T$-периодического решения
$\widetilde{x}_{\varepsilon},$ сходящегося при $\varepsilon\to 0$
к $\widetilde{x}^{\theta_0},$ является условие
\begin{equation}\label{suf_mi}
  M'(\theta_0)\not=0,
\end{equation}
то есть условие простоты корня $\theta_0.$

Различие между ситуациями, изучаемыми Малкиным и Мельниковым,
связаны с условиями, предъявляемыми к  $T$-периодической линейной
системе
\begin{equation}\label{ls_mi}
  \dot y=f'(\widetilde{x}(t))y.
\end{equation}

\noindent Так Малкин предполагал, что

{$$ \begin{array}{ll} \mbox{($C_{MA}$):} &  \mbox {алгебраическая
кратность мультипликатора +1 системы (\ref{ls_mi}) равна 1,}
\end{array}$$}

%{$$ \begin{array}{ll} \mbox{$(C_A)$:} & \mbox{цикл
%}\widetilde{x}\mbox{ изолирован и}
%\\
%& \mbox {алгебраическая кратность мультипликатора +1 системы
%(\ref{ls_mi}) равна 1,}
%\end{array}$$}

\noindent а Мельников, что

{$$ \begin{array}{ll} \mbox{$(C_{ME})$:} & \mbox {алгебраическая
кратность мультипликатора +1 системы (\ref{ls_mi}) равна 2,}
\\
& \mbox {но геометическая кратность мультипликатора +1 системы
(\ref{ls_mi}) равна 1.}
\end{array}$$}

%{$$ \begin{array}{ll} \mbox{$(C_{ME})$:} &  \mbox {геометическая
%кратность мультипликатора +1 системы (\ref{ls_mi}) равна 1.}
%\end{array}$$}

\noindent В статье будет показано, что другая идея Пуанкаре, об
индексе кривой относительно векторного поля, позволяет получить
близкие результаты при одном лишь предположении

{$$ \begin{array}{ll} \mbox{$(C)$:} & \mbox{в малой окрестности
цикла }\widetilde{x}\mbox{ нет других
}T\hskip-0.1cm-\hskip-0.2cm\mbox{ периодических решений}
\\
& \mbox {порождающей системы.}
\end{array}$$}

Условие ($C$)  выполнено как в ситуации Малкина ($C_{MA}$), так и
в ситуации Мельникова ($C_{ME}$). Но предположение ($C$) не
использует каких-либо свойств линеаризованной системы
(\ref{ls_mi}) и, следовательно, может быть выполнено для
вырожденных случаев, когда каждое решение системы (\ref{ls_mi})
является $T$-периодическим, см. \ref{degres}

Предлагаемый подход, в отличии от результатов Малкина и
Мельникова, не предполагает дифференцируемости возмущения и связан
с дальнейшим развитием теории индекса Пуанкаре, см. \cite{krazab}.
Хотя основная теорема доказывается для случая непрерывного
возмущения, следуя стандартным схемам (см., напр., \cite{book},
\S~6.2), она может быть перенесена не только на случай разрывных
возмущений, но и многозначных с выпуклыми образами.

 Будет показано, что условия существования периодических решений у возмущенной
системы вблизи цикла $\widetilde{x}$  могут быть связаны с числом
оборотов, которые делает вектор
\begin{equation}\label{Phi_mi}
  \Phi\left(\widetilde{x}(\theta)\right)=
M_E(\theta)\dot{\widetilde{ x}}(\theta)^\bot+
M_A(\theta)\dot{\widetilde{x}}(\theta),
\end{equation}
когда $\theta$ изменяется от $0$ до $T.$ Здесь использовано
обозначение
$$\xi^\bot= \left(
      \begin{array}{c}-\xi_2 \\ \xi_1
            \end{array}
      \right).$$  Это число оборотов называется индексом
Пуанкаре цикла $\widetilde{x}$ по отношению к полю $\Phi,$
заданному на цикле $\widetilde{x},$  и обозначается ${\rm
ind}(\widetilde{x},\Phi).$ Мы не останавливаемся здесь на
определении индекса Пуанкаре, которое имеется, например, в
(\cite{poin}, Гл.~III и Гл.~XIV), (\cite{lefs}, Гл.~IX, \S 4),
(\cite{andr}, Гл.~V, \S 10.4). Поскольку в указанной литературе
индекс Пуанкаре вводится для положительно ориентированной кривой,
то в дальнейшем мы всюду предполагаем, что кривая $\widetilde{x}$
является положительно ориентированной. Функции $M_A$ и $M_E$
формулы (\ref{Phi_mi}), которые определяются в следующем
параграфе, совпадают при некоторых дополнительных условиях с
классическими функциями Малкина и Мельникова.
%Отметим,
%наконец, что запись $"\theta\in[0,T]",$ использованная в формуле
%(\ref{Phi_mi}), означает, что равенство этой формулы справедливо
%при любом $\theta\in[0,T].$ Такой же смысл будет иметь указанное
%обозначение и в дальнейших формулах статьи.
\subsection{Основной результат}
\subsubsection{}\vskip-1cm
\addtocounter{sect}{1}

Для определения функций $M_E$ и $M_A$ рассмотрим  решения
следующей сопряженной системы
\begin{equation}\label{ss_mi}
  \dot z=-\left(f'(\widetilde{x}(t))\right)^*z.
\end{equation}
Именно, обозначим через $\widetilde{z}$ решение системы
(\ref{ss_mi}) с начальным условием $$\widetilde{z}(0)=\dot
{\widetilde{x}}(0)^\bot,$$ а через $\widehat {z}$ -- решение
системы (\ref{ss_mi}) с начальным условием
$$
  \widehat{z}(0)=\frac{1}{\|\dot{\widetilde{x}}(0)\|^2}\dot{\widetilde{x}}(0).
$$
Итак, положим
\begin{eqnarray}
  M_E^s(\theta)&=&\int\limits_{s-T+\theta}^{s+\theta}\left<\widetilde{z}(\tau),g(\tau-\theta,\widetilde{x}(\tau))\right>d\tau,\nonumber\\
  M_A^s(\theta)&=&\int\limits_{s-T+\theta}^{s+\theta}\left<\widehat{z}(\tau),g(\tau-\theta,\widetilde{x}(\tau))\right>d\tau, \qquad\theta,s\in[0,T],\nonumber
\end{eqnarray}
где $<\cdot,\cdot>$ -- обычное скалярное произведение в
$\mathbb{R}^2.$ При этом по определению считаем
$$M_E:=M_E^0,\quad M_A:=M_A^0.$$

\vskip0.5cm

Переходим к формулировке основной теоремы.

\begin{thm}\label{t1_mi}\hskip-0.2cm.
 Пусть выполнено  условие ($C$)  и возмущение в (\ref{ps_mi})
непрерывно. Предположим, что
\begin{equation}\label{L_cond}\tag{$A$}
  \mbox{для любых}\ s,\theta\in[0,T]\mbox{\ \ равенство\ \ } M_E^s(\theta)=0\mbox{\ \ влечет\ \
  }M_A^s(\theta)\not=0.
\end{equation}
Тогда, при достаточно малых $\varepsilon>0$ всякое
$T$-периодическое решение $\widetilde{x}_\varepsilon$ возмущенной
системы (\ref{ps_mi}) необходимо таково, что
\begin{equation}\label{us2_mi}
  \widetilde{x}_\varepsilon(t)\not=\widetilde{x}(s)\quad\mbox{при всех
  }t,s\in[0,T].
\end{equation}
Если же дополнительно имеем
\begin{equation}\label{ind_mi}\tag{B}
 {\rm ind}(\widetilde{x},\Phi)\not=1,
\end{equation}
то при всех достаточно малых $\varepsilon>0$ возмущенная система
(\ref{ps_mi}) действительно имеет, по крайней мере, два
$T$-периодических решения $\widetilde{x}_{\varepsilon,1}$ и
$\widetilde{x}_{\varepsilon,2},$ удовлетворяющих условию
(\ref{us2_mi}). Полученным решениям соответствуют
$\theta_1,\theta_2\in[0,T]$ такие, что
$$
  \widetilde{x}_{\varepsilon,1}(t)\to
  \widetilde{x}^{\theta_1}(t)\quad\mbox{\rm при}\ \varepsilon\to 0
$$
и
$$
  \widetilde{x}_{\varepsilon,2}(t)\to
  \widetilde{x}^{\theta_2}(t)\quad\mbox{\rm при}\ \varepsilon\to 0.
$$
Более того, решение $\widetilde{x}_{\varepsilon,1}$ лежит  внутри
цикла $\widetilde{x},$ а решение $\widetilde{x}_{\varepsilon,2}$
снаружи.
\end{thm}

Готовимся к доказательству теоремы~\ref{t1_mi}. Введем для этого
необходимые понятия и установим некоторые вспомогательные факты.

Обозначим через $\Omega(\cdot,t_0,\xi)$ решение $x$ порождающей
системы (\ref{np_mi}) с начальным условием $x(t_0)=\xi.$
Рассмотрим следующую систему (см. \cite{mal}, формула~3.8)
\begin{equation}\label{vsp_1}
  \dot y=f'_{(2)}(t,\Omega(t,0,\xi))y+g(t,\Omega(t,0,\xi))
\end{equation}
и обозначим через $\eta(\cdot,s,\xi)$ решение $y$ этой системы с
начальным условием $y(s)=0.$ Здесь и далее через $F'_{(k)}$
обозначается производная функции $F$ по $k$-й переменной. Поставим
в соответствие системе (\ref{ps_mi}) следующий оператор
$\widetilde{\Phi}^s:\mathbb{R}^2\to\mathbb{R}^2$ (см. \cite{dan},
теорема~1)
$$
  \widetilde{\Phi}^s(\xi)=\eta(T,s,\xi)-\eta(0,s,\xi).
$$
Отметим, что оператор $\widetilde{\Phi}^0$ не совпадает с $\Phi,$
однако мы покажем позже, что индексы Пуанкаре указанных
операторов, посчитанные на цикле $\widetilde{x},$ совпадают. В
случае, когда $f=0,$ то есть порождающей системой является система
$\dot x=0,$ имеем
$$
  \widetilde{\Phi}^s(\xi)=\int\limits_0^T g(\tau,\xi)d\tau,\qquad s\in\mathbb{R}, \xi\in\mathbb{R}^2.
$$
Следовательно, если $f=0,$ то $\widetilde{\Phi}^s$ является, с
точностью до множителя $1/T,$ классическим оператором усреднения
Крылова-Боголюбова-Митропольского возмущенной системы
(\ref{ps_mi}) (см. \cite{mit}, формула 7.112).

Обозначим через $\widehat{y}$ решение линеаризованной системы
(\ref{ls_mi}) с начальным условием
$$
  \widehat{y}(0)=\frac{1}{\|\dot{\widetilde{x}}(0)\|^2}\dot{\widetilde{x}}(0)^\bot
$$

Имеет место следующее свойство.

\begin{lem}\hskip-0.2cm. \label{Phi1_mi} Справедлива формула
\begin{equation}\label{Ftheta}
\widetilde{\Phi}^s(\widetilde
{x}(\theta))=M_A^s(\theta)\dot{\widetilde{ x}}(\theta)+
M_E^s(\theta)\widehat {y}(\theta), \qquad s,\theta\in\mathbb{R}.
\end{equation}
\end{lem}

Для доказательства леммы, а также для дальнейших рассуждений
используется утверждение, которое мы сейчас формулируем.

\begin{lem}\hskip-0.2cm.\label{form_eta} (см. \cite{non}, лемма~1 или \cite{can}, лемма~2) Имеет место формула
$$ \eta(\theta,s,\xi)=\Omega'_{(3)}(\theta,0,\xi)\int\limits_s^\theta
\Omega'_{(3)}(0,t,\Omega(t,0,\xi))g(t,\Omega(t,0,\xi)) d\tau, $$
соответственно
$$
\widetilde{\Phi}^s(\xi)=\int\limits_{s-T}^s
\Omega'_{(3)}(0,\tau,\Omega(\tau,0,\xi))g(\tau,\Omega(\tau,0,\xi))
d\tau.
$$
\end{lem}

\noindent {\rm Д о к а з а т е л ь с т в о\ \ \ л
 е м м ы\ \ \ \ref{Phi1_mi}.}\ \  В силу леммы~\ref{form_eta} имеем
$$
\widetilde{\Phi}^s(\xi)=\int\limits_{s-T}^s
  (\Omega)'_{(3)}(0,\tau,\Omega(\tau,0,\xi))g(\tau,\Omega(\tau,0,\xi))d\tau.
$$
Далее (см., напр., \cite{kraope}, теорема~2.1),
$\Omega'_{(3)}(t,0,x_0(\theta))=Y(t,\theta),$ где $Y(t,\theta)$ --
фундаментальная матрица системы
\begin{equation}\label{tp}
\dot y(t) = f'(\widetilde{x}(t+\theta))y(t),
\end{equation}
удовлетворяющая условию $Y(0,\theta)=I$ и так как
$\Omega'_{(3)}(0,t,\Omega(t,0,\xi))\Omega'_{(3)}(t,0,\xi)=I$ для
любых $t\in\mathbb{R},$ $\xi\in\mathbb{R}^2,$ имеем
\begin{equation}\label{good}
  \widetilde{\Phi}^s(\widetilde{x}(\theta))=
  \int\limits_{s-T}^s
  Y^{-1}(\tau,\theta)g(\tau,\widetilde{x}(\tau+\theta))d\tau,\qquad
  s,\theta\in[0,T].
\end{equation}
Покажем теперь, что
\begin{equation}\label{cl1}
  Y^{-1}(\tau,\theta)=Y(\theta,0) Y^{-1}(\tau+\theta,0),\qquad
  \tau,\theta\in[0,T].
\end{equation}
Действительно, легко видеть, что $Y(t+\theta,0)$ --
фундаментальная матрица системы (\ref{tp}) и, таким образом,
$Y(t+\theta,0)Y^{-1}(\theta,0)$ также является фундаментальной
матрицей для (\ref{tp}), более того имеем
$Y(t+\theta,0)Y^{-1}(\theta,0)=I$ в $t=0.$ Следовательно,
$Y(t+\theta,0)Y^{-1}(\theta,0)=Y(t,\theta),$ что эквивалентно
(\ref{cl1}).

\noindent Подставляя (\ref{cl1}) в (\ref{good}) и используя замену
переменных $\tau+\theta=t$ в интеграле выражения (\ref{good}),
получаем
\begin{eqnarray}
\widetilde{\Phi}^s(\widetilde{x}(\theta))& = &
\int\limits_{s-T}^s
  Y^{-1}(\tau,\theta)g(\tau,\widetilde{x}(\tau+\theta))d\tau=
  Y(\theta,0)\int\limits_{s-T}^{s}
  Y^{-1}(\tau+\theta,0)g(\tau,\widetilde{x}(\tau+\theta))d\tau=\nonumber\\
& =&  Y(\theta,0)\int\limits_{s-T+\theta}^{s+\theta}
  Y^{-1}(t,0)g(t-\theta,\widetilde{x}(t))dt.\nonumber
\end{eqnarray}
Положим $\widetilde {Z}(t)=(\widehat {z}(t),\widetilde {z}(t))$ и
обозначим через $Z(t)$ фундаментальную матрицу сопряженной системы
(\ref{ss_mi}) такую, что $Z(0)=I,$ имеем $Z(t)=\widetilde
{Z}(t)\widetilde {Z}^{-1}(0).$ Здесь и далее использована запись
$(a,b),$ где $a,b\in\mathbb{R}^2,$ для обозначения матрицы, первым
столбцом которой является $a$ и вторым $b.$ По лемме Перрона (см.
\cite{per} или \cite[Sec. III, \S 12]{dem}) $Y^{-1}(t)=Z^*(t)$ и
значит
\begin{eqnarray}
\Phi^s(\widetilde {x}(\theta))&=&
Y(\theta)\int\limits_{s-T+\theta}^{s+\theta}
  Y^{-1}(\tau)g(\tau-\theta,\widetilde{
  x}(\tau))d\tau=\nonumber\\
& =& (\widetilde{
Z}^*(\theta))^{-1}\int\limits_{s-T+\theta}^{s+\theta}
  \widetilde {Z}^*(\tau)g(\tau-\theta,\widetilde{
  x}(\tau))d\tau= \nonumber\\
&=&  (\widetilde{ Z}^*(\theta))^{-1}\left(\begin{array}{c}
M_A^s(\theta)\\
  M_E^s(\theta)
\end{array}\right).\nonumber
\end{eqnarray}
%(\widetilde {x}(t)\ \widetilde {y}(t))
Но, замечая, что
$$  (\widetilde{
Z}^*(0))^{-1}=\left(\left(\begin{array}{ll}
   \left(1/{\|\dot{\widetilde{x}}(0)\|^2}\right)\dot{\widetilde{x}}_1(0)
   & -\dot{\widetilde{x}}_2(0) \\
   \left(1/{\|\dot{\widetilde{x}}(0)\|^2}\right)\dot{\widetilde{x}}_2(0) & \dot{\widetilde{x}}_1(0)
  \end{array}\right)^{-1}\right)^*=$$
$$=
  \frac{1}{{\rm det}\left\|\left(\begin{array}{ll}
   \left(1/{\|\dot{\widetilde{x}}(0)\|^2}\right)\dot{\widetilde{x}}_1(0)
   & -\dot{\widetilde{x}}_2(0) \\
   \left(1/{\|\dot{\widetilde{x}}(0)\|^2}\right)\dot{\widetilde{x}}_2(0) & \dot{\widetilde{x}}_1(0)
  \end{array}\right)
  \right\|}
  \circ $$
  $$
  \qquad\qquad\circ\left(\begin{array}{ll}
   \dot{\widetilde{x}}_1(0)
   & \dot{\widetilde{x}}_2(0) \\
   -\left(1/{\|\dot{\widetilde{x}}(0)\|^2}\right)\dot{\widetilde{x}}_2(0) &
   \left(1/{\|\dot{\widetilde{x}}(0)\|^2}\right)\dot{\widetilde{x}}_1(0)
  \end{array}\right)^*=$$
  $$=\left(\begin{array}{ll}
   \dot{\widetilde{x}}_1(0)
   & -\left(1/{\|\dot{\widetilde{x}}(0)\|^2}\right)\dot{\widetilde{x}}_2(0) \\
   \dot{\widetilde{x}}_2(0) &
   \left(1/{\|\dot{\widetilde{x}}(0)\|^2}\right)\dot{\widetilde{x}}_1(0)
  \end{array}\right)=\left(\dot{\widetilde{x}}(0),{\widehat{y}}(0)\right),
$$
заключаем
$$
  (\widetilde{
  Z}^*(\theta))^{-1}=\left(\dot{\widetilde{x}}(\theta),{\widehat{y}}(\theta)\right)\quad
  \mbox{для любого }\theta\in\mathbb{R}.
$$

Лемма доказана.

В качестве следующего шага на пути к доказательству
теоремы~\ref{t1_mi} мы устанавливаем, что ${\rm
ind}(\widetilde{x},\widetilde{\Phi})={\rm
ind}(\widetilde{x},\Phi).$ Как и в случае функций $M^s_E$ и
$M^s_A$ считаем по определению
$$
  \widetilde{\Phi}:=\widetilde{\Phi}^0.
$$

\begin{lem}\label{defdef}\hskip-0.2cm. Пусть выполнено условие (\ref{L_cond}) теоремы~\ref{t1_mi}.
Тогда
$$
  {\rm ind}(\widetilde{x},\Phi)={\rm ind}(\widetilde{x},\widetilde{\Phi}).
$$
\end{lem}

\noindent {\rm Д о к а з а т е л ь с т в о. } Определим деформацию
\begin{eqnarray}
 \Phi_\lambda(\widetilde{x}(\theta))&=&\left(\int\limits_{-T+\lambda\theta}^{\lambda\theta}
\left<\widehat{z} (\tau),g(\tau-\theta,\widetilde{
  x}(\tau))\right>d\tau\right)\dot{\widetilde{
  x}}(\theta)+\nonumber\\
  & &+
  \left(\int\limits_{-T+\lambda\theta}^{\lambda\theta}
\left<\widetilde{z} (\tau),g(\tau-\theta,\widetilde{
  x}(\tau))\right>d\tau\right)
\left(\lambda\widehat{y}(\theta)+(1-\lambda)\dot{\widetilde{x}}(\theta)^\bot\right),\quad
\theta\in[0,T],\nonumber
\end{eqnarray} Так как
$\Phi_1(\widetilde{x}(\theta))=\widetilde{\Phi}(\widetilde{x}(\theta))$
и $\Phi_0(\widetilde{x}(\theta))={\Phi}(\widetilde{x}(\theta)),$
то достаточно установить, что
\begin{equation}\label{diff}
  \Phi_\lambda(\widetilde{x}(\theta))\not=0\quad\mbox{для всех
  }\lambda\in[0,1],\ \theta\in[0,T].
\end{equation}
По определению $\widehat{y}$ имеем
$\left<\widehat{y}(0),\dot{\widetilde{x}}(0)^\bot\right>>0.$ Но в
силу линейной независимости векторов $\dot{\widetilde{x}}(\theta)$
и $\widehat{y}(\theta)$ для любого $\theta\in[0,T]$ также имеем
$\left<\widehat{y}(\theta),\dot{\widetilde{x}}(\theta)^\bot\right>\not=0.$
Поэтому
$$
  \left<\widehat{y}(\theta),\dot{\widetilde{x}}(\theta)^\bot\right>>0\quad\mbox{для
  любого }\theta\in[0,T]
$$
и, значит,
\begin{equation}\label{zna}
  \left<\lambda\widehat{y}(\theta)+(1-\lambda)\dot{\widetilde{x}}(\theta)^\bot,
  \dot{\widetilde{x}}(\theta)^\bot\right>>0\quad\mbox{для
  любого }\theta\in[0,T],\ \lambda\in[0,T].
\end{equation}
Предположим, что (\ref{diff}) не верно, то есть существуют
$\lambda_0\in[0,1],$ $\theta_0\in[0,T]$ такие, что
$$
\Phi_{\lambda_0}(\widetilde{x}(\theta_0))=0.
$$
Учитывая  (\ref{zna}), из последнего равенства заключаем, что
$$
 \int\limits_{-T+\lambda_0\theta_0}^{\lambda_0\theta_0}
\left<\widehat{z} (\tau),g(\tau-\theta_0,\widetilde{
  x}(\tau))\right>d\tau=0\quad\mbox{и}\quad \int\limits_{-T+\lambda_0\theta_0}^{\lambda_0\theta_0}
\left<\widetilde{z} (\tau),g(\tau-\theta_0,\widetilde{
  x}(\tau))\right>d\tau=0.
$$
Но $\lambda_0\theta_0\in[0,T]$ и, следовательно, последнее
соотношение противоречит предположению (\ref{L_cond}).

Лемма доказана.

Наконец, нам необходимо следующее утверждение.

\begin{lem}\hskip-0.2cm.\label{zamena}
Функция $x\in C([0,T],\mathbb{R}^2)$ является $T$-периодическим
решением системы (\ref{ps_mi}) тогда и только тогда, когда функция
\begin{equation}\label{zp_}
  \nu(t)=\Omega(0,t,x(t)),\quad t \in \mathbb{R}
\end{equation}
является решением системы
$$
\dot\nu=\varepsilon
\Omega'_{(3)}(0,t,\Omega(t,0,\nu))g(t,\Omega(t,0,\nu)),\quad
t\in\mathbb{R},
$$
удовлетворяющим условию $\nu(0)=\Omega(T,0,\nu(T)).$
\end{lem}

\noindent {\rm Д о к а з а т е л ь с т в о. } Произведем в системе
(\ref{ps_mi}) замену переменных
\begin{equation}\label{zpn}
  x(t)=\Omega(t,0,\nu(t)).
\end{equation}
Формула  (\ref{zpn}) каждому $\nu \in C([0,T],\mathbb{R}^2)$
ставит в соответствие $x \in C([0,T],\mathbb{ R}^2) $ гомеоморфно,
и обратное отображение дается формулой (\ref{zp_}). Следовательно,
функция $x$ является решением системы (\ref{ps_mi}) тогда и только
тогда, когда функция $\nu,$ введенная по закону (\ref{zp_}),
удовлетворяет следующему равенству
\begin{eqnarray}\label{f9}
  \Omega'_{(1)}(t,0,\nu(t))+\Omega'_{(3)}(t,0,\nu(t))\dot \nu(t)=\varepsilon g(t,\Omega(t,0,\nu(t)))+ \nonumber \\
  +f(\Omega(t,0,\nu(t))),\qquad t\in\mathbb{R}.
\end{eqnarray}
По определению функции $\Omega$ имеем
$$
\Omega'_{(1)}(t,0,\nu(t))=f(\Omega(t,0,\nu(t))),
$$
поэтому, учитывая, что
$$
  \Omega'_{(3)}(0,t,\Omega(t,0,\xi))\Omega'_{(3)}(t,0,\xi)=I,\qquad
  t\in\mathbb{R}, \ \xi\in\mathbb{R}^2,
$$
 система (\ref{f9}) может быть переписана в виде
\begin{equation}\label{1_}
  \dot \nu(t)=\varepsilon
  \Omega'_{(3)}(0,t,\Omega(t,0,\nu(t)))g(t,\Omega(t,0,\nu(t))),\qquad t\in\mathbb{R}.
\end{equation}
Рассмотрим произвольное $T$-периодическое решение $x$ системы
(\ref{ps_mi}). Имеем
$$
  \nu(0)=\Omega(0,0,x(0))=x(0)=x(T)=\Omega(T,0,\nu(T)).
$$

Лемма доказана.

Обозначим через $U\subset\mathbb{R}^2$ внутренность цикла
$\widetilde{x}.$

Для доказательства теоремы~\ref{t1_mi} рассмотрим вполне
непрерывный интегральный оператор
$Q_\varepsilon:C([0,T],\mathbb{R}^2)\to C([0,T],\mathbb{R}^2)$
\begin{equation}\label{operator}
(Q_\varepsilon x)(t)=x(T)+\int\limits_0^t
f(x(\tau))d\tau+\varepsilon\int\limits_0^t g(\tau,x(\tau))d\tau,
\quad t\in[0,T]
\end{equation}
на множестве
\[
W_U=\left\{\widehat {x}\in C([0,T],\mathbb{R}^2):\widehat{
x}(t)\in U,\ \mbox{\rm для\ любого\ } t\in [0,T]\right\}.
\]
Через $d(I-F,W)$ будем обозначать
 степень Лерэ-Шаудера преобразования $I-F:\overline{W}\to\mathbb{R}^2$  относительно нуля (см. \cite{leray}, Гл. I, \S 5). Иногда, чтобы подчеркнуть, что $F$ задано в пространстве
 $E$ мы будем писать $d_E(I-F,W).$

\

Основную роль в доказательстве теоремы~\ref{t1_mi} играет
следующая теорема.

\begin{thm}\hskip-0.2cm.\label{fromnon}
Если

 $\widetilde{\Phi}^s(\xi)\not=0$ для любого $\xi\in\partial U$ и любого
$s\in[0,T],$

\noindent то существует $\varepsilon_0>0$ такое, что при
$\varepsilon\in(0,\varepsilon_0]$ справедливы следующие
утверждения

1) для любого $x\in C([0,T],\mathbb{R}^2)$ такого, что
$Q_\varepsilon x=x$ имеем $x(t)\not\in\partial U$ для всех
$t\in[0,T],$ в частности, оператор $Q_\varepsilon$ не имеет
неподвижных точек на $\partial W_U;$

2) $
  d(I-Q_\varepsilon,W_U)={\rm ind}(\widetilde{x},\widetilde{\Phi}),\quad
  \varepsilon\in(0,\varepsilon_0].
$
\end{thm}

Утверждение 2) теоремы~\ref{fromnon} впервые предложено в
дипломной работе \cite{diplom} автора и опубликовано в \cite{non}
(теорема~2). Формулировка результата о существовании периодических
решений для возмущенной системы (\ref{ps_mi}), непосредственно
следующего из утверждения 2), впервые опубликована в
\cite{dan}(теорема~1). Поскольку доказательство одного лишь
утверждения 1) почти совпадает с доказательством обоих утверждений
1) и 2), нам показалось целесообразным привести ниже полное
доказательство теоремы~\ref{fromnon}. Тем более, оно значительно
проще, чем приведенное в \cite{non}, где возмущенная система имеет
два содержащих $\varepsilon>0$ слагаемых и $\varepsilon>0$ входит
в эти слагаемые с разными степенями.

Теорема~\ref{fromnon} является развитием результатов Дж.~Мавена,
который установил (см. \cite{maw1} и \cite{maw2}) похожие
утверждения в случае нулевой или линейной порождающей системы.
Хотя Мавен рассматривал случай пространства $\mathbb{R}^n,$
доказательство теоремы~\ref{fromnon} немедленно переносится и на
$n$-мерный случай.

\noindent {\rm Д о к а з а т е л ь с т в о\ \ \ т
 е о р е м ы\ \ \ \ref{fromnon}.} \ \ Положим
$$\Upsilon(t,\xi)=
\Omega'_{(3)}(0,t,\Omega(t,0,\xi))g(t,\Omega(t,0,\xi)).
$$
Согласно лемме~\ref{zamena} каждой неподвижной точке из $W_U$
оператора $Q_\varepsilon$ соответствует неподвижная точка
(\ref{zp_}) оператора
$$
(G_\varepsilon
\nu)(t)=\Omega(T,0,\nu(T))+\int\limits_0^t\Upsilon(\tau,\nu(\tau))d\tau,
$$
которая, как легко проверить, вновь принадлежит $W_U.$
Следовательно, если $d(I-G_\varepsilon,W_U)$ определен, то, в силу
(\cite{krazab}, теорема 26.4), имеем
$$
  d(I-Q_\varepsilon,W_U)=d(I-G_\varepsilon,W_U).
$$
В пространстве $ C([0,T],\mathbb{ R}^2)$ рассмотрим
вспомогательный вполне непрерывный оператор
$$
  (A_{\varepsilon}\nu)(t)=\Omega(T,0,\nu(T))-
    \varepsilon\int \limits_0^T\Upsilon(\tau,\nu(\tau))d\tau
$$
и покажем, что при достаточно малых $\varepsilon>0$ поля
$I-G_\varepsilon$ и $I-A_{\varepsilon}$ гомотопны на границе
множества $W_U$. Зададим следующую деформацию
$$
D_{\varepsilon}(\lambda,\nu)(t)=\nu(t)-\Omega(T,0,\nu(T))-\varepsilon\int
\limits_0^{\lambda t+(1-\lambda)T}\Upsilon(\tau,\nu(\tau))d\tau,
$$
соединяющую поля $G_\varepsilon$ и $G_{1,\varepsilon}$. Покажем,
что при достаточно малых $\varepsilon>0$ деформация
$D_\varepsilon$ невырождена на границе множества $W_U.$ Для этого
будет доказано более сильное утверждение, которое будет
использовано затем для доказательства утверждения 1), а именно
покажем, что существует $\varepsilon_0>0$ такое, что при всех
$\varepsilon\in(0,\varepsilon_0]$ и $\lambda\in[0,1]$ каждое
решение уравнения $D_\varepsilon(\lambda,\nu)=\nu$ удовлетворяет
условию $\nu(t)\not\in\partial U$ при всех $t\in[0,T].$
Предположим, что это не так. Тогда для произвольной
последовательности
${\left\{\varepsilon_k\right\}}_{k=1}^{\infty}\subset(0,1]$ такой,
что $\varepsilon\to 0$ при $k\to\infty$ найдутся
последовательности
${\left\{\lambda_k\right\}}_{k=1}^{\infty}\subset[0,1]$ и
${\left\{\nu_k\right\}}_{k=1}^{\infty}\subset C([0,T],\mathbb{
R}^2),$ при которых
\begin{equation}\label{f15}
  \nu_k(t)=\Omega(T,0,\nu_k(T))+\varepsilon_k\int \limits_0^{\lambda_k t+(1-\lambda_k)T}
  \Upsilon(\tau,\nu_k(\tau))d\tau,\ t\in[0,T]
\end{equation}
и
\begin{equation}\label{nuzhno}
  \nu_k([0,T])\cap\partial U\not=\emptyset.
\end{equation}
 Так как последовательность чисел
${\left\{\lambda_k\right\}}_{k=1}^{\infty}$ ограничена, то из нее
можно выделить сходящуюся подпоследовательность. Поэтому,  без
ограничения общности можем считать, что последовательность
${\left\{\lambda_k\right\}}_{k=1}^{\infty}$ сходится. Из
(\ref{nuzhno}) следует, что функции последовательности
${\left\{\nu_k\right\}}_{k=1}^{\infty}$  равномерно ограничены.
Поэтому, на основании непрерывности функции $\Upsilon$ найдется
константа $c_0>0$ такая, что
$\left\|\Upsilon(t,\nu_k(t))\right\|\le c_0$, $t\in[0,T],$
$k\in\mathbb{ N},$ и для любых $t_1,\ t_2\in[0,T],$
$k\in\mathbb{N}$ имеем оценку
$$
  \|\nu_k(t_2)-\nu_k(t_1)\|=\varepsilon_k\left\|\int
    \limits_{\lambda_k t_1+(1-\lambda_k)T}^{\lambda_k t_2+
     (1-\lambda_k T)}\Upsilon(\tau,z_k(\tau))d\tau\right\|\le
       \varepsilon_k\lambda_k(t_2-t_1)c_0,
$$
из которой следует, что функции последовательности
${\left\{\nu_k\right\}}_{k=1}^{\infty}$ равностепенно непрерывны.
Значит, применяя теорему Арцела, из этой последовательности можно
выделить сходящуюся подпоследовательность. Поэтому мы без
ограничения общности можем считать, что последовательность
${\left\{\nu_k\right\}}_{k=1}^{\infty}$ сходится. Положим
$\lambda_0=\lim_{k\to\infty}\lambda_k$ и
$\nu_0=\lim_{k\to\infty}\nu_k.$ Тогда $\lambda_0\in[0,1]$ и
$\nu_0([0,T])\cap\partial U\not=\emptyset.$ Так как $\dot \nu_k\to
0$ при $k\to\infty,$ то функция $\nu_0$ постоянна. Соотношение
(\ref{nuzhno}) эквивалентно существованию числа $t_k\in[0,T]$
такого, что $\nu_k(t_k)\in\partial U.$ Тогда
\begin{equation}\label{f21}
 \Omega(T,0,\nu_k(t_k))=\nu_k(t_k), \quad k\in\mathbb{N}.
\end{equation}
Вычитая из равенства (\ref {f15}), записанного при $t=T$, это же
равенство, записанное при $t=t_n$, получаем
\begin{equation}\label{f22}
  \nu_k(T)-\nu_k(t_k)=\varepsilon_k\int \limits_{\lambda_k t_k+(1-\lambda_k)T}^T\Upsilon(\tau,\nu_k(\tau))d\tau.
\end{equation}
На основании (\ref{f21}) выражение (\ref{f15}) при $t=T$ можно
переписать в виде
\begin{eqnarray}
  \nu_k(T)-\nu_k(t_k)=\Omega(T,0,\nu_k(T))-\Omega(T,0,\nu_k(t_k))+\nonumber\\
  +\varepsilon_k\int \limits_0^T\Upsilon(\tau,z_k(\tau))d\tau,\nonumber
\end{eqnarray}
или
\begin{eqnarray}\label{f23}
  (I-\Omega'_{(3)}(T,0,\nu_k(t_k)))(\nu_k(T)-\nu_k(t_k))= \nonumber\\
  =\varepsilon_k\int \limits_0^T\Upsilon(\tau,\nu_k(\tau))d\tau+o(\nu_k(t_k),\nu_k(T)-\nu_k(t_k)),
\end{eqnarray}
где функция $o(\xi,h)$ удовлетворяет соотношению
\begin{equation}
  \frac{\|o(\xi,h)\|}{\|h\|} \to 0,\ \mbox{\rm при}\ \|h\| \to 0,\ \xi\in \mathbb{R}^2.
\end{equation}
Подставляя (\ref {f22}) в (\ref {f23}), получаем равенство
\begin{eqnarray}\label{f25}
  (I-\Omega'_{(3)}(T,0,\nu_k(t_k)))\int \limits_{\lambda_k t_k+(1-\lambda_k)T}^T\Upsilon(\tau,\nu_k(\tau))d\tau =
  \nonumber \\
  =\int \limits_0^T\Upsilon(\tau,\nu_k(\tau))d\tau+\frac{o(\nu_k(t_k),\nu_k(T)-\nu_k(t_k))}{\varepsilon_k}.
\end{eqnarray}
Из (\ref {f22}) следует, что найдется константа $c>0$ такая, что
$$
\|\nu_k(T)-\nu_k(t_k)\| \le c\varepsilon_k.
$$
Откуда
\begin{equation}\label{f26}
  \left\|\frac{o(\nu_k(t_k),\nu_k(T)-\nu_k(t_k))}{\varepsilon_k} \right\| \le c \frac{\|o(\nu_k(t_k),\nu_k(T)-\nu_k(t_k))\|}
  {\|\nu_k(T)-\nu_k(t_k) \|}.
\end{equation}
Из (\ref {f26}), учитывая то, что значения функций $\nu_k$
ограничены равномерно по $k\in\mathbb{ N},$ следует, что
\begin{equation}\label{f27}
  \left\|\frac{o(\nu_k(t_k),\nu_k(T)-\nu_k(t_k))}{\varepsilon_k} \right\| \to 0 \ \mbox{\rm при}\ k \to \infty.
\end{equation}
Совершив, учитывая (\ref{f27}), предельный переход при $k \to
\infty$ в (\ref {f25}), получим
$$
  (I-\Omega'_{(3)}(T,0,\xi_0))\int \limits_s^T\Upsilon(\tau,\xi_0)d\tau =
  \int \limits_0^T\Upsilon(\tau,\xi_0)d\tau,
$$
или
\begin{eqnarray}\label{f28}
  \int\limits_s^T\Omega'_{(3)}(T,0,\xi_0)\Omega'_{(3)}(0,\tau,\Omega(\tau,0,\xi_0))g(\tau,\Omega(\tau,0,\xi_0))d\tau- \nonumber\\
  -\int\limits_s^0\Omega'_{(3)}(0,\tau,\Omega(\tau,0,\xi_0))g(\tau,\Omega(\tau,0,\xi_0))d\tau=0,
\end{eqnarray}
где $s=\lim_{k\to\infty}\left(\lambda_k
t_k+(1-\lambda_k)T\right)\in[0,T].$ Пользуясь
леммой~\ref{form_eta} равенство (\ref{f28}) можно переписать в
виде
$$
  \eta(T,s,\xi_0)-\eta(0,s,\xi_0)=0,
$$
в чем противоречие с предположением теоремы. Таким образом,
существует $\varepsilon_0
> 0$ такое, что при всех $\varepsilon\in(0,\varepsilon_0]$ и
$\lambda\in[0,1]$ каждое решение уравнения
$D_\varepsilon(\lambda,\nu)=\nu$ удовлетворяет условию
$\nu(t)\not\in\partial U$ при всех $t\in[0,T].$ При $\lambda=1$
полученный результат совпадает с утверждением 1) теоремы. Перейдем
к доказательству утверждения 2). Как уже говорилось, доказанное
свойство означает, в частности, что
\begin{equation}\label{f31}
  D_\varepsilon(\lambda,\nu)\not= 0,\  \lambda \in [0,1] , \ \nu \in \partial W_U, \
  \varepsilon\in(0,\varepsilon_0],
\end{equation}
то есть поля $I-G_{\varepsilon}$ и $I-A_{\varepsilon}$ гомотопны
на границе множества $W_U$ при $\varepsilon\in(0,\varepsilon_0].$
Обозначим через $C_0([0,T],\mathbb{R}^2)$ подпространство
пространства $C([0,T],\mathbb{ R}^2),$ состоящее из всех
постоянных функций, определенных на отрезке $[0,T]$ и принимающих
значения в $\mathbb{R}^2.$ Имеем $A_{\varepsilon}(\partial
W_U)\subset C_0([0,T],\mathbb{ R}^2).$ Далее, по построению
множество $W_U$ содержит функции, тождественно равные
произвольному фиксированному элементу из $U.$ Наконец, из
(\ref{f31}) при $\lambda=0$ получаем
$$
  A_{\varepsilon}(\nu)\not= \nu,\  \nu \in \partial W_U, \ \varepsilon\in(0,\varepsilon_0],
$$
откуда, учитывая соотношение $\partial (W_U\cap
C_0([0,T],\mathbb{R}^2))\subset\partial W_U,$ следует, что при
$\varepsilon\in(0,\varepsilon_0]$ поле $I-A_{\varepsilon}$ не
имеет нулей на границе множества $W_U\cap C_0([0,T],\mathbb{
R}^2).$ Поэтому, при $\varepsilon\in(0,\varepsilon_0]$ законно
сужение поля $I-A_{\varepsilon}$ на подпространство
$C_0([0,T],\mathbb{ R}^2),$ что означает
\begin{equation} \label{vr1}
   d_{C([0,T],{\mathbb{R}}^2)}(I-A_{\varepsilon},W_U)=
     d_{C_0([0,T],\mathbb{ R}^2)}(I-A_{\varepsilon},
       W_U\cap C_0([0,T],\mathbb{ R}^2)),
\end{equation}
где в левой и правой частях равенства записаны топологические
степени в пространствах $C([0,T],\mathbb{ R}^2)$ и
$C_0([0,T],\mathbb{ R}^2)$ соответственно.

Заметим, что постоянная функция $\nu\in W_U\cap C_0([0,T],\mathbb{
R}^2)$ тогда и только тогда является решением уравнения
$A_{\varepsilon}\nu=\nu,$ когда элемент $\xi=\nu(0)$ является
решением уравнения $A^0_{\varepsilon}\xi=\xi,$ где
$$
  A_{\varepsilon}^0\xi=\Omega(T,0,\xi)+\varepsilon\int
  \limits_0^{T}\Upsilon(\tau,\xi)d\tau.
$$
Применяя теорему об эквивалентных уравнениях к уравнениям
$A_{\varepsilon}\nu=\nu$ и $A_{\varepsilon}^0\xi=\xi,$ а затем
теорему о переходе к подпространству, при
$\varepsilon\in(0,\varepsilon_0]$ получаем
\begin{equation}\label{vr2}
    { d}_{C_0([0,T],\mathbb{R}^2)}(I-A_{\varepsilon},W_U\cap C_0([0,T],\mathbb{R}^2))
     = {d}_{\mathbb{R}^2}(I-A_\varepsilon^0,U).
\end{equation}
Для вычисления топологической степени ${d}_{\mathbb{
R}^2}(I-A_\varepsilon^0,U)$ положим
$$
  A_{1,\varepsilon}\xi=-\varepsilon\int \limits_0^{T}\Upsilon(\tau,\xi)d\tau.
$$
Так как $ (I-A_{\varepsilon}^0)(\xi)=A_{1,\varepsilon}\xi,\ $
$\xi\in\partial U,$ то при $\varepsilon\in(0,\varepsilon_0]$ имеем
\begin{equation}\label{vr3}
  {d}_{\mathbb{R}^2}(I-A_{\varepsilon}^0,U)=
    {d}_{\mathbb{ R}^2}(A_{1,\varepsilon},U).
\end{equation}
Покажем, что векторные поля $A_{1,\varepsilon}$ и $A_{1,1}$
гомотопны на границе множества $U$ при
$\varepsilon\in(0,\varepsilon_0].$ Зададим линейную деформацию
$$
  D_{1,\varepsilon}(\lambda,\xi)=(\lambda\varepsilon+1-\lambda)\int \limits_0^{T}\Upsilon(\tau,\xi)d\tau,\ \ \xi\in U
$$
и установим, что она невырождена на границе множества $U.$
Предположим противное, тогда для некоторых $\lambda\in[0,1],\ $
$\xi\in\partial U$ и $\varepsilon\in(0,\varepsilon_0)$ будем иметь
$$
  (\lambda_0\varepsilon+1-\lambda)\int \limits_0^{T}\Upsilon(\tau,\xi)d\tau=0,
$$
откуда
$$
  \int \limits_0^{T}\Upsilon(\tau,\xi)d\tau=0,
$$
что, в силу леммы~\ref{form_eta}, противоречит условиям теоремы.
Таким образом, пользуясь леммой~\ref{form_eta},
\begin{equation}\label{vr4}
  { d}_{\mathbb{ R}^2}(A_{1,\varepsilon},U)=
     {d}_{\mathbb{ R}^2}(A_{1,1},U)={ d}_{\mathbb{ R}^2}(-\eta(T,0,\cdot),U).
\end{equation}
Подставляя (\ref{vr4}) в (\ref{vr3}), получаем
\begin{equation}\label{res_}
  { d}_{\mathbb{ R}^2}(I-A_\varepsilon^0,U)=
    {d}_{\mathbb{ R}^2}(-\eta(T,0,\cdot),U),\ \varepsilon\in(0,\varepsilon_0).
\end{equation}
Подставляя (\ref{vr2}) в (\ref{vr1}), пользуясь гомотопностью
полей $I-A_\varepsilon$ и $I-G_\varepsilon$ и соотношением
(\ref{res_}), окончательно имеем
$$
  { d}_{C([0,T],\mathbb{ R}^2)}(I-G_\varepsilon,W_U)=
      { d}_{\mathbb{ R}^2}(-\eta(T,0,\cdot),U), \ \varepsilon\in(0,\varepsilon_0].
$$
Но поле $\eta(T,0,\cdot)$ получается из поля $-\eta(T,0,\cdot)$
непрерывным поворотом всех векторов на $180^0$ против часовой
стрелки, следовательно
$$
{ d}_{\mathbb{ R}^2}(-\eta(T,0,\cdot),U)={ d}_{\mathbb{
R}^2}(\eta(T,0,\cdot),U).
$$
 По определению множества $U$ и, пользуясь
независимостью топологической степени от способа ее определения
(см. \cite{krazab}, утверждение с.~17), имеем ${ d}_{\mathbb{
R}^2}(\eta(T,0,\cdot),U)={\rm
ind}(\widetilde{x},\eta(T,0,\cdot))={\rm
ind}(\widetilde{x},\widetilde{\Phi}).$

Теорема доказана.

Все готово для того, чтобы перейти к доказательству
теоремы~\ref{t1_mi}.

\noindent {\rm Д о к а з а т е л ь с т в о\ \ \ т
 е о р е м ы\ \ \ \ref{t1_mi}.} \ \ Пусть $\varepsilon_0>0$ то, о котором
говорится в теореме \ref{fromnon}, тогда, учитывая также
леммы~\ref{Phi1_mi} и \ref{defdef},
$$
  d(Q_\varepsilon,W_U)={\rm ind}(\widetilde{x},\Phi)\quad \mbox{для любых
  }\varepsilon\in(0,\varepsilon_0]
$$
и пользуясь предположением (\ref{ind_mi}),
\begin{equation}\label{dee1}
  d(Q_\varepsilon,W_U)\not=1,\quad \mbox{для любых
  }\varepsilon\in(0,\varepsilon_0].
\end{equation}
 Положим $U_\delta^-=U\backslash B_\delta(\partial
U),$ $U_\delta^+=U\cup B_\delta(\partial U).$ Здесь и далее через
$B_\delta(D)$ обозначена $\delta$-окрестность множества $D$ в
норме содержащего $D$ пространства. На основании условия ($C$)
можно зафиксировать такое $\delta_0>0,$ что порождающая система
(\ref{np_mi}) не имеет $T$-периодических решений с начальными
условиями из $\partial U_\delta^-\cup\partial U_\delta^+$ при всех
$\delta\in(0,\delta_0].$ Без ограничения общности можем считать,
что $\delta_0>0$ выбрано достаточно малым так, что
$U_\delta^-\not=\emptyset.$ По теореме Капетто-Мавена-Занолина
(\cite{camaza}, Следствие~1) имеем
$$
   d(Q_0,W_{U_\delta^-})=d_{\mathbb{R}^2}(f,U_\delta^-)\mbox{ и }
   d(Q_0,W_{U_\delta^+})=d_{\mathbb{R}^2}(f,U_\delta^+)\quad\mbox{при
   всех }\delta\in(0,\delta_0].
$$
Без ограничения общности можно считать, что малость $\delta_0>0$
достаточна для того, чтобы
$$
  d_{\mathbb{R}^2}(f,U_\delta^-)=d_{\mathbb{R}^2}(f,U_\delta^+)=d_{\mathbb{R}^2}(f,U),\quad
  \delta\in(0,\delta_0].
$$
По теореме Пуанкаре (см. Лефшец \cite[теорема 11.1]{lefs} или
Красносельский и др.  \cite[теорема~2.3]{kradr}) имеем
$d_{\mathbb{R}^2}(f,U)=1,$ поэтому
$$
   d(Q_0,W_{U_\delta^-})=
   d(Q_0,W_{U_\delta^+})=1\quad\mbox{при
   всех }\delta\in(0,\delta_0].
$$
Таким образом, каждому $\delta\in(0,\delta_0]$ соответствует
$\varepsilon_\delta>0$ такое, что
\begin{equation}\label{dee2}
   d(Q_\varepsilon,W_{U_\delta^-})=
   d(Q_\varepsilon,W_{U_\delta^+})=1\quad\mbox{при
   всех }\varepsilon\in(0,\varepsilon_\delta] \mbox{ и }\delta\in(0,\delta_0].
\end{equation}
Без ограничения общности можно считать, что
$\varepsilon_\delta<\varepsilon_0$ при всех
$\delta\in(0,\delta_0].$ Тогда, из (\ref{dee1}) и (\ref{dee2})
получаем, что при всех $\delta\in(0,\delta_0]$ и
$\varepsilon\in(0,\varepsilon_\delta]$ система (\ref{ps_mi})
имеет, по крайней мере, два $T$-периодических решения
$x_{\varepsilon,1}\in W_U\backslash W_{U_\delta^-}$ и
$x_{\varepsilon,2}\in W_{U_\delta^+}\backslash W_U.$ Из этого, в
частности, имеем $x_{\varepsilon,1}(0)\in U,$
$x_{\varepsilon,2}(0)\not\in U$ и, используя утверждение 1)
теоремы \ref{fromnon}, заключаем, что $x_{\varepsilon,1}(t)\in U,$
$x_{\varepsilon,2}(t)\not\in U$ при всех $t\in[0,T].$

 Теорема доказана.

\subsection{Вычисление индекса}
\subsubsection{}\vskip-1cm
\addtocounter{sect}{1}

В настоящем параграфе показывается, что если функция $M$ имеет
ровно два нуля на интервале $[0,T),$ то проверка условия
(\ref{ind_mi}) отличия от $+1$ индекса Пуанкаре, участвующего в
формулировке теоремы~\ref{t1_mi}, сводится к проверке
алгебраического неравенства. Основным будет являться следующее
утверждение.

\begin{thm}\hskip-0.2cm. \label{analt} Пусть возмущение в (\ref{ps_mi})
непрерывно.
 Предположим, что
функция $M_E$ имеет ровно два нуля $\theta_1$ и $\theta_2$ на
интервале $[0,T).$ Тогда, если $M_E$ строго монотонна в точках
$\theta_1$ и $\theta_2$ и
$$
  M_A(\theta_1)\cdot M_A(\theta_2)<0,
$$
то либо ${\rm ind}(\widetilde{x},\Phi)=0,$ либо ${\rm
ind}(\widetilde{x},\Phi)=2.$

\end{thm}

Доказательство теоремы основано на следующей лемме.

\begin{lem}\hskip-0.2cm. (см. \cite{anal}, Лемма~5) \label{anal} Пусть
$q:[0,T]\to\mathbb{R}^2,$ $q(0)=q(T)$ -- жорданова кривая и
$\psi:\mathbb{R}^2\to\mathbb{R}^2$ -- непрерывное векторное поле
такое, что $\psi(q(t))\not=0$ для каждого $t\in[0,T].$
Предположим, что существует направляющая функция
$z:[0,T]\to\mathbb{R}^2,$ $z(0)=z(T)$ такая, что:

1) $\left<z(\theta),\dot {q}(\theta)\right>\not=0$ для каждого
$\theta\in[0,T],$

2) скалярная функция
$\alpha(\theta)=\left<\psi(q(\theta)),z(\theta)\right>$ имеет
ровно два нуля $\theta_1,\theta_2$ на интервале $[0,T)$ и строго
монотонна в этих точках,

3) ${\rm sign}\left<\psi(q(\theta_1)), \left(\begin{array}{c}
z_2(\theta_1) \\ -z_1(\theta_1)
\end{array}\right)
\right>=-{\rm sign}\left<\psi(q(\theta_2)),\left(\begin{array}{c}
z_2(\theta_2) \\ -z_1(\theta_2)
\end{array}\right)\right>.$

Тогда либо ${\rm ind}(q,\psi)=0,$ либо ${\rm ind}(q,\psi)=2.$

\end{lem}

Отметим, что лемма~\ref{anal} является уточнением теоремы
Борсука-Улама  \cite{borsuk} о четных векторных полях для
рассматриваемого специального случая.

\vskip0.3cm

\noindent {\rm Д о к а з а т е л ь с т в о\ \ \ т
 е о р е м ы\ \ \ \ref{analt}.} \ \ Используем лемму~\ref{anal}. Для этого
выберем $q=\widetilde{x},$ $\psi=\Phi,$ $z=\dot{\widetilde{x}}.$
Имеем
$$\alpha(\theta)=-M_E(\theta)\quad\mbox{и}\quad
\left<\Phi(\widetilde{x}(\theta)), \left(\begin{array}{c} \dot{\widetilde{x}}_2(\theta) \\
-\dot{\widetilde{x}}_1(\theta)
\end{array}\right)
\right> =
-\left<\Phi(\widetilde{x}(\theta)),\dot{\widetilde{x}}(\theta)^\bot\right>=-M_A(\theta).$$
То есть условия леммы~\ref{anal} совпадают с предположениями
доказываемой теоремы, что завершает доказательство.

\begin{rem} Верна теорема, полученная из теоремы~\ref{analt} заменой $M_A$
на $M_E$ и, соответственно, $M_E$ на $M_A.$
\end{rem}

\begin{rem}
Аналогичные теореме~\ref{analt} утверждения могут быть доказаны в
случае, когда функция $M_E$ имеет произвольное число нулей. В этом
случае необходимо требовать, чтобы знаки функции $M_A$ в этих
нулях подходящим образом чередовались.
\end{rem}

\subsection{Сопоставление с теоремами Малкина и Мельникова}
\subsubsection{}\vskip-1cm
\addtocounter{sect}{1}

В первую очередь мы установим, что функция $M_A$ совпадает с
функцией Малкина в случае, когда выполнены условия Малкина
($C_{MA}$), а функция $M_E$ с функцией Мельникова в случае, когда
выполнены условия Мельникова ($C_{ME}$) и некоторые дополнительные
свойства симметрии.

 В работе \cite{mal}
И.~Г.~Малкин в предположении ($C_{MA}$) определяет бифуркационную
функцию $\widetilde{M}_A$ как (см. \cite{mal}, формула~3.13)
$$
  \widetilde{M}_A(\theta)=\int\limits_0^{T}\left<\widetilde{\widetilde{z}}(\tau),g(\tau-\theta,\widetilde{x}(\tau))\right>d\tau,
$$
где $\widetilde{\widetilde{z}}$ -- $T$-периодическое решение
сопряженной системы (\ref{ss_mi}), удовлетворяющее условию
\begin{equation}\label{normi}
  \left<\widetilde{\widetilde{z}}(0),\dot{\widetilde{x}}(0)\right>=1.
\end{equation}
В \cite{mal} отмечается, что такой выбор всегда возможен (см.
\cite{mal}, формулы~2.13 и 3.7). При условии ($C_{MA}$) такой
выбор необходимо единственен, то есть функция Малкина определена
однозначно. В этой же работе Малкин предлагает следующий результат
(см. \cite{mal}, утверждение~с.~638).

\vskip0.3cm

 \noindent{\rm Т е о р е м а  \ \ М а л к и н а.}\ \
{\it Пусть правая часть возмущенной системы (\ref{ps_mi})
непрерывно дифференцируема и выполнены условия Малкина ($C_{MA}$).
Если возмущенная система (\ref{ps_mi}) имеет при достаточно малых
$\varepsilon>0$ $T$-периодическое решение
$\widetilde{x}_{\varepsilon}$ такое, что
\begin{equation}\label{convma}
  \widetilde{x}_\varepsilon(t)\to
  \widetilde{x}^{\theta_0}(t)\quad\mbox{\rm при}\ \varepsilon\to
  0,
\end{equation}
то
 $$\widetilde{M}_A(\theta_0)=0.$$ Если же вместе с указанным необходимым
условием нуль $\theta_0$ является простым, то есть
${\widetilde{M}_A}{}'(\theta_0)\not=0,$ то при всех достаточно
малых $\varepsilon>0$ возмущенная система (\ref{ps_mi})
действительно имеет $T$-периодическое решение
$\widetilde{x}_\varepsilon,$ удовлетворяющее (\ref{convma}). }

\vskip0.3cm

Следующее утверждение дает условие, при которых $M_A^s(\theta)$ не
зависит от $s$ и совпадает с ${\widetilde{M}}_A.$

\begin{lem}\label{MAlem} Пусть выполнено условие Малкина ($C_{MA}$).
Если $\widehat{y}$ является собственной функцией линеаризованной
системы (\ref{ls_mi}), то решение $\widehat{z}$ --
$T$-периодическое. В частности,
$$
  M_A^s(\theta)={\widetilde{M}}_A(\theta)\quad\mbox{при всех
  }\ s,\theta\in[0,T].
$$
\end{lem}

\noindent {\rm Д о к а з а т е л ь с т в о.} Пусть
$\widetilde{\widetilde{z}}$ -- $T$-периодическая функция,
участвующая в определении функции $\widetilde{M}_A,$ и,
следовательно, удовлетворяющая (\ref{normi}). Так как
$\widehat{y}$ -- собственная функция системы (\ref{ls_mi}), то,
учитывая условие ($C_{MA}$), существует $\rho\not=1$ такое, что
$$
  \widehat{y}(T)=\rho\widehat{y}(0).
$$
Но, в силу леммы Перрона (см. \cite{per} или \cite[Sec. III, \S
12]{dem}),
$$
  \left<\widetilde{\widetilde{z}}(0),\widehat{y}(0)\right>=
  \left<\widetilde{\widetilde{z}}(T),\widehat{y}(T)\right>,
$$
что возможно только в случае, когда
$$
  \left<\widetilde{\widetilde{z}}(0),\widehat{y}(0)\right>=0.
$$
С другой стороны по определению $\widehat{z}$
$$
  \left<\widehat{z}(0),\dot{\widetilde{x}}(0)\right>=1\quad\mbox{и}\quad
  \left<\widehat{z}(0),\widehat{y}(0)\right>=0.
$$
Таким образом, векторы $\widehat{z}(0)$ и
$\widetilde{\widetilde{z}}(0)$ в скалярном произведении как с
вектором $\dot{\widetilde{x}}(0),$ так и с вектором
$\widehat{y}(0)$ принимают одни и те же значения. Следовательно
$\widehat{z}(0)=\widetilde{\widetilde{z}}(0).$

Лемма доказана.

\vskip0.3cm

Таким образом, при выполнении условий леммы~\ref{MAlem}  проверка
условий теоремы~\ref{t1_mi} упрощается.

%\begin{cor} Eсли выполнено  условие $(C_{MA})$  и  $\widehat{y}$ является собственной функцией
%линеаризованной системы (\ref{ls_mi}), то для справедливости
%условия (\ref{L_cond}) достаточно предположить выполнение
%следующего свойства
%\begin{equation}\tag{$A_{MA}$}
%\mbox{для любого }\theta_0\in[0,T] \mbox{ такого, что
%}\widetilde{M}_A(\theta_0)=0, \mbox{ имеем }
%M_E^s(\theta_0)\not=0\mbox{ при всех}\ s\in[0,T].
%\end{equation}
%\end{cor}

Обратимся теперь к случаю Мельникова ($C_{ME}$). В работе
\cite{mel} Мельников вводит функцию (см. \cite{mel}, формула для
$A_0(v),$ c.~42)
$$
\widetilde{M}_E(\theta)= \int\limits_0^{T}{\rm
det}\left\|\left(\dot{\widetilde{x}}(\tau),g(\tau-\theta,\widetilde{x}(\tau))\right)\right\|d\tau.
$$
и устанавливает следующий результат (см. \cite{mel}, Лемма~7).

\vskip0.3cm

{\rm Т е о р е м а\ \  М е л ь н и к о в а.} \ \ {\it Пусть правая
часть возмущенной системы (\ref{ps_mi}) дважды непрерывно
дифференцируема и цикл $\widetilde{x}$ вложен в некоторое
семейство циклов порождающей гамильтоновой системы (\ref{np_mi}).
Пусть выполнены условия Мельникова ($C_{ME}$). Если возмущенная
система (\ref{ps_mi}) имеет при достаточно малых $\varepsilon>0$
$T$-периодическое решение $\widetilde{x}_{\varepsilon},$
сходящееся к циклу $\widetilde{x}$ при $\varepsilon\to 0$ с
точностью до сдвижки, то существует $\theta_0\in[0,T]$ такое, что
 $$
 \widetilde{M}_E(\theta_0)=0.
 $$
 Если же вместе с указанным необходимым
условием нуль $\theta_0$ является простым, то есть
${\widetilde{M}_E}{}'(\theta_0)\not=0,$ то при всех достаточно
малых $\varepsilon>0$ возмущенная система (\ref{ps_mi})
действительно имеет $T$-периодическое решение
$\widetilde{x}_\varepsilon$ такое, что выполнено (\ref{convma}). }

На самом деле, в цитированной работе Мельников рассматривал случай
аналитических правых частей, но для интересующего нас утверждения
действительно достаточно требовать двойной непрерывной
дифференцируемости, соответствующая формулировка и доказательство
имеется, например, в \cite{guck}, теорема~4.6.2.

Для функции $M_E^s$ условия ее независимости от $s$ и совпадения с
$\widetilde{M}_E$ получаются неодинаковыми. К выводу указанных
условий мы переходим.

\vskip0.3cm

 Установим вначале одно вспомогательное утверждение.

\begin{lem}\hskip-0.2cm.\label{lemma1}
Предположим, что $T$-периодическая система
\begin{equation}\label{t1}
  \dot u =A(t)u,\quad u\in\mathbb{R}^2
\end{equation}
имеет мультипликатор $+1$ алгебраической кратности 2, и
$\widetilde {u}$ -- $T$-периодическое решение этой системы такое,
что
$$ \widetilde {u}_1(0)=0,\ \widetilde {u}_2(0)\not=0.$$
 Тогда для решения $\widehat {u}$ системы
(\ref{t1}), удовлетворяющего условию
$$\widehat {u}_1(0)\not=0,\ \widehat
{u}_2(0)=0,$$ справедлива формула
$$
\widehat {u}(t+T)=\widehat {u}(t)+\frac{\widehat
{u}_2(T)}{\widetilde {u}_2(0)}\widetilde {u}(t),\quad
t\in\mathbb{R}.
$$
\end{lem}

\noindent {\rm Д о к а з а т е л ь с т в о.}
 Обозначим через $X$ нормированную ($X(0)=I$) фундаментальную
матрицу системы (\ref{t1}). Так как
$X(T_0)\left(\begin{array}{c}0\\1\end{array}\right)=
\left(\begin{array}{c}0\\1\end{array}\right),$ то $X(T_0)=
\left(\begin{array}{cc}a&0\\b&1\end{array}\right).$ По условию
леммы $X(T)$ имеет два собственных значения $+1,$ значит $X(T_0)=
\left(\begin{array}{cc}1&0\\b&1\end{array}\right),$ где
$b\in\mathbb{R}$ -- некоторое число. Имеем
\begin{eqnarray}
  X(t+T_0)\widehat {u}(0)&=&X(t)X(T_0)\widehat{
  u}(0)=X(t)\left(\begin{array}{cc}1&0\\b&1\end{array}\right)\widehat{
  u}(0)=X(t)\widehat {u}(0)+X(t)\left(\begin{array}{c}0\\b\widehat{
  u}_1(0)\end{array}\right)=\nonumber\\
  &=&X(t)\widehat {u}(0)+\frac{b\widehat {u}_1(0)}{\widetilde {u}_2(0)}
  \widetilde{
  u}(t).\nonumber
\end{eqnarray}
В то же время
$$
  X(T_0)\widehat{
  u}(0)=\left(\begin{array}{cc}1&0\\b&1\end{array}\right)\widehat{
  u}(0)=\widehat {u}(0)+\left(\begin{array}{c}0\\b\widehat{
  u}_1(0)\end{array}\right),
$$
откуда $b\widehat {u}_1(0)=\widehat {u}_2(T).$ Лемма доказана.

\vskip0.3cm

Следующая лемма утверждает, что выполнение условия Мельникова
($C_{ME}$) достаточно для независимости $M_E^s(\theta)$  от $s$ (в
случае Малкина были нужны  дополнительные условия, см.
лемму~\ref{lemma1}).

\begin{lem}\hskip-0.2cm. \label{lemma2} Если алгебраическая кратность мультипликатора
$+1$ линеаризованной системы (\ref{ls_mi}) равна $2,$  то решение
$\widetilde {z}$ является $T$-периодическим. В частности,
$M_E^s(\theta)$ не зависит от $s.$
\end{lem}

Условие леммы~\ref{lemma2} выполнено как в случае Мельникова
($C_{ME}$), так и в вырожденном случае, когда каждое решение
системы (\ref{ls_mi}) является $T$-периодическим.

 \noindent {\rm Д
о к а з а т е л ь с т в о\ \ \ л е м м ы}\ \ \ \ref{lemma2}. \ \ \
Пусть $\theta\in[0,T]$ таково, что
\begin{equation}\label{takovo}
  \dot{\widetilde{x}}_1(\theta)=0.
\end{equation}
Тогда функция $\dot{\widetilde{x}}{}^{\theta}$ является решением
системы
\begin{equation}\label{lss_mi}
  \dot y=f'\left(\widetilde{x}(t+\theta)\right)y
\end{equation}
с начальными условием $\dot{\widetilde{x}}{}^\theta_1(0)=0.$
Обозначим через  $\widehat{\widehat{y}}$ решение системы
(\ref{lss_mi}) с начальным условием
$\widehat{\widehat{y}}(0)=(1,0).$ Используя  лемму~\ref{lemma1},
заключаем, что
\begin{equation}\label{zakl}
  \widehat{\widehat{y}}(T)=\widehat{\widehat{y}}(0)+\frac{\widehat{\widehat{y}}_2(T)}
  {\dot{\widetilde{x}}_2{}^\theta(0)}\dot{\widetilde{x}}{}^\theta(T).
\end{equation}
Если $\widehat {\widehat{y}}_2(T)=0,$ то из (\ref{zakl}) имеем,
что каждое решение системы (\ref{lss_mi}), а значит и системы
(\ref{ss_mi}), является $T$-периодическим.
 Рассмотрим случай, когда
\begin{equation}\label{ca}
\widehat {\widehat{y}}_2(T)\not=0.
\end{equation} В силу теоремы о периодических
решениях сопряженной системы (см. \cite[Гл. III, \S 23, теорема
2]{dem}), система
\begin{equation}\label{sss_mi}
  \dot z=-\left(f'\left(\widetilde{x}(t+\theta)\right)\right)^*z
\end{equation}
имеет по крайней мере одно $T$-периодическое решение, обозначим
это решение через $\widetilde{\widetilde{z}}.$ Домножая равенство
(\ref{zakl}) скалярно на $\widetilde{\widetilde{z}},$ получаем
$$
    \left<\widehat{y}^\theta(T),\widetilde{\widetilde{z}}(T)\right>=
    \left<\widehat{y}^\theta(0),\widetilde{\widetilde{z}}(T)\right>+\frac{\widehat{y}_2^\theta(T)}
  {\dot{\widetilde{x}}_2{}^\theta(0)}\left<\dot{\widetilde{x}}{}^\theta(t),\widetilde{\widetilde{z}}(T)\right>
  =
  \left<\widehat{y}^\theta(0),\widetilde{\widetilde{z}}(0)\right>+\frac{\widehat{y}_2^\theta(T)}
  {\dot{\widetilde{x}}_2{}^\theta(0)}\dot{\widetilde{x}}{}^\theta_2(0)\widetilde{\widetilde{z}}_2(0).
$$
Но в силу леммы Перрона (см. \cite{per} или \cite[Sec. III, \S
12]{dem}) $
\left<\widehat{y}^\theta(T),\widetilde{\widetilde{z}}(T)\right>=
\left<\widehat{y}^\theta(0),\widetilde{\widetilde{z}}(0)\right>,$
поэтому $\widetilde{\widetilde{
 z}}_2(0)=0.$ Учитывая (\ref{takovo}), заключаем
$
  \left<\dot{\widetilde{x}}{}^\theta(0),\widetilde{\widetilde{z}}(0)\right>=0
$ или
$$
  \left<\dot{\widetilde{x}}(0),\widetilde{\widetilde{z}}{}^{-\theta}(0)\right>=0.
$$
 Следовательно, векторы $\widetilde{\widetilde{
 z}}{}^{-\theta}(0)$ и $\widetilde{z}(0)$ линейно зависимы, то
 есть существует $a\not=0$ такое, что $\widetilde{z}(0)=a \widetilde{\widetilde{
 z}}{}^{-\theta}(0).$
 Но обе функции $\widetilde{\widetilde{
 z}}{}^{-\theta}$ и $\widetilde{z}$
 является решениями одной и той же сопряженной системы
 (\ref{ss_mi}), поэтому
 $$\widetilde{z}(t)=a\widetilde{\widetilde{
 z}}{}^{-\theta}(t)\quad\mbox{при всех }t\in[0,T],$$
в частности, функция $\widetilde{z}$ -- $T$-периодическая.

Лемма доказана.

\vskip0.3cm
 Наконец, предъявим условия, при которых
$M_E=\widetilde{M}_E.$

\begin{lem}\hskip-0.2cm.\label{symlem} Пусть выполнено условие Мельникова
($C_{ME}$). Предположим, что для всех
  $\xi\in\mathbb{R}^2$ выполнено
\begin{equation}\label{m1}
  f_1(\xi)=f_1(-\xi_1,\xi_2),
\end{equation}
\begin{equation}\label{m2}
  f_2(\xi)=-f_2(-\xi_1,\xi_2),
\end{equation}
\begin{equation}\label{m3}
  (f_1)'_{(1)}(\xi)=-(f_2)'_{(2)}(\xi).
\end{equation}
Тогда для любого $\theta\in\mathbb{R}$ имеем
$$
  \widehat{z}(\theta)=\left(\begin{array}{c} \widehat {y}_2(\theta)
 \\ -\widehat
 {y}_1(\theta)\end{array}\right), \qquad \widetilde{z}(\theta)=\left(\begin{array}{c} -\dot
{\widetilde {x}}_2(\theta)
 \\\dot {\widetilde {x}}_1(\theta)\end{array}\right)
$$
и, в частности,
\begin{eqnarray}
  M_E^s(\theta)&=&\widetilde{M}_E(\theta),\nonumber\\
  M_A^s(\theta)&=&-\int\limits_{s-T}^{s}{\rm
det}\left\|\left({\widehat{y}}(\tau),g(\tau-\theta,\widetilde{x}(\tau))\right)\right\|d\tau,
\qquad s,\theta\in[0,T].\nonumber
\end{eqnarray}
\end{lem}
Для доказательства леммы~\ref{symlem} будет использовано следующее
утверждение.

\begin{lem}\hskip-0.2cm. \label{lemma7} Рассмотрим линейную систему
\begin{equation}\label{lin1}
  \left(\begin{array}{l}
    \dot y_1\\
    \dot y_2
  \end{array}\right)=\left(\begin{array}{cc}
    a(t) & d(t)\\
    b(t) & -a(t)
  \end{array}\right)\left(\begin{array}{l}
    y_1\\
    y_2
  \end{array}\right)
\end{equation}
Предположим, что $a(-t)=-a(t),$ $b(-t)=b(t),$ $d(-t)=d(t)$ для
любых $t\in[c_1,c_2].$ Тогда, если $y$ -- некоторое решение
системы (\ref{lin1}), то функция $z(t)=(y_2(-t),y_1(-t))$
удовлетворяет на отрезке $[c_1,c_2]$ сопряженной к (\ref{lin1})
системе.
\end{lem}

Справедливость леммы~\ref{lemma7} проверяется непосредственно
подстановкой решения $z(t)=(y_2(-t),y_1(-t))$ в сопряженную
систему.

\noindent {\rm Д о к а з а т е л ь с т в о\ \ \ л
 е м м ы\ \ \ \ref{symlem}.} \ \ Пусть $\tau\in[0,T]$ таково, что
$$
  \widetilde{x}_1^\tau(0)=0.
$$

 Пользуясь условиями (\ref{m1}) и (\ref{m2}), легко
проверить, что функция
$p(t)=(-\widetilde{x}_1^\tau(-t),\widetilde{x}_2^\tau(-t))$
является решением порождающей системы (\ref{np_mi}). Но
$p(0)=\widetilde{x}^\tau(0),$ следовательно
\begin{equation}\label{sym1}
 (-\widetilde{x}_1^\tau(-t),\widetilde{x}_2^\tau(-t))=\widetilde{x}^\tau(t) \mbox{ для любого
 }t\in[0,T].
\end{equation}
Рассмотрим линеаризованную на $\widetilde{x}^\tau$ порождающую
систему (\ref{np_mi})
\begin{equation}\label{slin_mi}
  \left(\begin{array}{l}
    \dot y_1\\
    \dot y_2
  \end{array}\right)=\left(\begin{array}{cc}
    (f_1)'_{(1)}(\widetilde{x}^\tau(t)) & (f_1)'_{(2)}(\widetilde{x}^\tau(t))\\
    (f_2)'_{(1)}(\widetilde{x}^\tau(t)) & (f_2)'_{(2)}(\widetilde{x}^\tau(t))
  \end{array}\right)\left(\begin{array}{l}
    y_1\\
    y_2
  \end{array}\right)
\end{equation}
Из условия (\ref{m3}) следует, что
$$
  (f_1)'_{(1)}(\widetilde{x}^\tau(t))=-(f_2)'_{(2)}(\widetilde{x}^\tau(t)).
$$
Из (\ref{m1}) имеем
$-(f_1)'_{(1)}(-\xi_1,\xi_2)=(f_1)'_{(1)}(\xi_1,\xi_2)$ и,
учитывая (\ref{sym1}), получаем
\begin{equation}\label{m5}
  (f_1)'_{(1)}(\widetilde{x}^\tau(t))=-(f_1)'_{(1)}(\widetilde{x}^\tau(-t)).
\end{equation}
Из (\ref{m2}) имеем
$-(f_2)'_{(1)}(\xi_1,\xi_2)=(f_2)'_{(1)}(-\xi_1,\xi_2)$ и,
учитывая (\ref{sym1}), получаем
\begin{equation}\label{m4}
  (f_2)'_{(1)}(\widetilde{x}^\tau(-t))=(f_2)'_{(1)}(\widetilde{x}^\tau(t)).
\end{equation}
Наконец, из (\ref{m1}) имеем
$(f_1)'_{(2)}(-\xi_1,\xi_2)=(f_1)'_{(2)}(\xi_1,\xi_2)$ и, учитывая
(\ref{sym1}), получаем
\begin{equation}\label{m6}
  (f_1)'_{(2)}(\widetilde{x}^\tau(t))=(f_1)'_{(2)}(\widetilde{x}^\tau(-t)).
\end{equation}
Таким образом, выполнены условия леммы~\ref{lemma7}, на основании
которой получаем, что функции
$$
 \widehat {\widehat{z}}{}^\tau(t)=\left(\begin{array}{c} \widehat {y}_2^\tau(-t)
 \\ \widehat {y}_1^\tau(-t)\end{array}\right)\quad\mbox{и}\quad
\widetilde {\widetilde{z}}{}^\tau(t)=\left(\begin{array}{c} \dot
{\widetilde {x}}_2{}^\tau(-t)
 \\\dot {\widetilde {x}}_1{}^\tau(-t)\end{array}\right)
$$
удовлетворяют сопряженной к (\ref{slin_mi}) системе
\begin{equation}\label{ssss}
  \dot z=-\left(f'(\widetilde{x}^\tau(t))\right)z.
\end{equation}
%и, вместе с тем, условию (\ref{initz}).
Из (\ref{sym1}) для любого $t\in[0,T]$ имеем
\begin{equation}\label{nzvezda}
\dot {\widetilde{x}}_1{}^\tau(-t)=\dot
{\widetilde{x}}_1{}^\tau(t),\quad -\dot
{\widetilde{x}}_2{}^\tau(-t)=\dot{\widetilde{x}}_2{}^\tau(t).
\end{equation}
Покажем, что вместе с $\widehat{y}^\tau$ решением
системы (\ref{slin_mi}) является функция
$p(t)=(-\widehat{y}_1^\tau(-t),\widehat{y}_2^\tau(-t)).$
Действительно, из (\ref{m5}) и (\ref{m6}) имеем
$$
 \dot
 p_1(t)=-(f_1)'_{(1)}(\widetilde{x}^\tau(t))y_1^\tau(-t)+(f_1)'_{(2)}(\widetilde{x}^\tau(t))y_2^\tau(-t),
$$
и из (\ref{m4}), (\ref{m3}) и (\ref{m5}) имеем
$$
 \dot
 p_2(t)=-(f_2)'_{(1)}(\widetilde{x}^\tau(t))y_1^\tau(-t)+(f_2)'_{(2)}(\widetilde{x}^\tau(t))y_2^\tau(-t).
$$
Вместе с тем, из того, что
$\left<\widehat{\widehat{z}}{}^0(0),\widetilde{x}^0(0)\right>=0$ в
силу леммы Перрона (см. \cite{per} или \cite[Sec. III, \S
12]{dem}) имеем
$\left<\widehat{\widehat{z}}{}^\tau(0),\widetilde{x}^\tau(0)\right>=0,$
то есть
$$\left<\left(\begin{array}{c}
\widehat{y}^\tau_2(0) \\
\widehat{y}^\tau_1(0)
\end{array}\right),\widetilde{x}^\tau(0)\right>=0.$$
Но $\widetilde{x}_1^\tau(0)=0,$ поэтому из последнего равенства
заключаем, что $\widehat{y}^\tau_1(0)=0.$ Полученное свойство
позволяет утверждать, что
 $p(0)=\widehat{y}^\tau(0),$ поэтому
\begin{equation}\label{symy}
  (-\widehat{y}_1^\tau(-t),\widehat{y}_2^\tau(-t))=\widehat{y}^\tau(t), \quad
  t\in\mathbb{R}.
\end{equation}
На основании (\ref{nzvezda}) и (\ref{symy}) функции
$\widehat{\widehat{z}}{}^\tau$ и
$\widetilde{\widetilde{z}}{}^\tau$ можно переписать в виде
$$
 \widehat {\widehat{z}}{}^\tau(t)=\left(\begin{array}{c} \widehat {y}_2^\tau(t)
 \\ -\widehat {y}_1^\tau(t)\end{array}\right)\quad\mbox{и}\quad
\widetilde {\widetilde{z}}{}^\tau(t)=\left(\begin{array}{c} -\dot
{\widetilde {x}}_2{}^\tau(t)
 \\\dot {\widetilde {x}}_1{}^\tau(t)\end{array}\right).
$$
Так как функции $\widehat{\widehat{z}}{}^\tau$ и
$\widetilde{\widetilde{z}}{}^\tau$ являются решениями сопряженной
системы (\ref{ssss}), то функции $\widehat{\widehat{z}}$ и
$\widetilde{\widetilde{z}}$ являются решениями сопряженной системы
(\ref{ss_mi}). Но
$$
  \widehat {\widehat{z}}(0)=\left(\begin{array}{c} \widehat {y}_2(0)
 \\ -\widehat
 {y}_1(0)\end{array}\right)=\frac{1}{\|\dot{\widetilde{x}}(0)\|^2}\dot{\widetilde{x}}(0)=\widehat{z}(0),
$$
$$
  \widetilde {\widetilde{z}}(0)=\left(\begin{array}{c} -\dot
{\widetilde {x}}_2(0)
 \\\dot {\widetilde {x}}_1(0)\end{array}\right)=\widetilde{z}(0),
$$
следовательно,
$$
  \widehat{z}(t)=\left(\begin{array}{c} \widehat {y}_2(t)
 \\ -\widehat
 {y}_1(t)\end{array}\right), \qquad \widetilde{z}(t)=\left(\begin{array}{c} -\dot
{\widetilde {x}}_2(t)
 \\\dot {\widetilde {x}}_1(t)\end{array}\right).
$$

Лемма доказана.

Таким образом, выполнение условий симметрии леммы~\ref{symlem} не
только приводит функцию $M_E$ к классической, но и упрощает
вычисление функции $M_A^s.$ Приведем одно следствие из
леммы~\ref{lemma1}, которое также может быть использовано для
упрощения вычисления функции $M_A^s.$

\begin{cor}\hskip-0.2cm. Пусть алгебраическая кратность мультипликатора
$+1$ линеаризованной системы (\ref{ls_mi}) равна $2$  и
$\widehat{z}_2(0)=0.$ Тогда
$$
  M_A^s(\theta)=M_A^T(0)-\frac{\widehat{z}_2(T)}{\widetilde{z}_2(0)}\int\limits_{s+\theta}^T
  \left<\widetilde{z}(\tau),g(\tau-\theta,\widetilde{x}(\tau))\right>d\tau.
$$
\end{cor}
\noindent {\rm Д о к а з а т е л ь с т в о. } Используя замену
переменных $t=\tau+T$ в интеграле и лемму~\ref{lemma1}, можем
провести следующее преобразование
$$
M_A^s(\theta)=\int\limits_{s-T+\theta}^{s+\theta}
  \left<\widehat {z}(\tau),g(\tau-\theta,\widetilde {x}(\tau))\right>d\tau=
$$
$$=
  \int\limits_0^{s+\theta}\left<\widehat {z}(\tau),g(\tau-\theta,\widetilde{
  x}(\tau))\right>d\tau+\int\limits_{s-T+\theta}^0\left<\widehat{
  z}(\tau),g(\tau-\theta,\widetilde{
  x}(\tau))\right>d\tau=
$$
$$
= \int\limits_0^{s+\theta}\left<\widehat{
z}(\tau),g(\tau-\theta,\widetilde{
  x}(\tau))\right>d\tau+\int\limits_{s+\theta}^{T}\left<\widehat{
  z}(t-T),g(t-\theta,\widetilde{
  x}(t))\right>dt=
$$
$$
= \int\limits_0^{s+\theta}\left<\widehat{
z}(\tau),g(\tau-\theta,\widetilde{
  x}(\tau))\right>d\tau+\int\limits_{s+\theta}^{T}\left<\left(\widehat {z}(t)-
  \frac{\widehat {z}_2(T)}{\widetilde {z}_2(0)}\widetilde {z}(t)\right),g(t-\theta,\widetilde
  {x}(t))\right>dt=
$$
$$
= \int\limits_0^{T}\left<\widehat
{z}(\tau),g(\tau-\theta,\widetilde{
  x}(\tau))\right>d\tau-
  \frac{\widehat {z}_2(T)}{\widetilde {z}_2(0)}\int\limits_{s+\theta}^{T}\left<\widetilde
  {z}(t),g(t-\theta,\widetilde{
  x}(t))\right>dt=
$$
$$
  =M_A^T(\theta)-\frac{\widehat {z}_2(T)}{\widetilde {z}_2(0)}\int\limits_{s+\theta}^{T}\left<\widetilde
  {z}(t),g(t-\theta,\widetilde{
  x}(t))\right>dt.
$$

Следствие доказано.

\vskip0.3cm

Закончив сопоставление функций $M_A^s$ и $M_E^s$ с классическими
функциями $\widetilde{M}_A$ и $\widetilde{M}_E$ соответственно
Малкина и Мельникова, переходим к сопоставлению
теоремы~\ref{t1_mi} с соответствующими классическими теоремами.
Именно, при помощи нескольких примеров будут сопоставлены
заключения, к каким приводит теорема Мельникова и к каким
теорема~\ref{t1_mi}.
 Совершенно понятно, что заменой порождающих
систем предлагаемых примеров на подходящие системы, допускающие
$T$-периодический предельный цикл, сопоставление
теоремы~\ref{t1_mi} может быть проведено и с теоремой Малкина. В
качестве такой подходящей системы может быть рассмотрена,
например, система
$$
  \left(\begin{array}{c}\dot x_1 \\ \dot x_2 \end{array}\right)=
  \left(\begin{array}{c} x_2+x_1(x_1^2+x_2^2-1) \\
    -x_1+x_2(x_1^2+x_2^2-1) \end{array}\right),
$$
допускающая единственный цикл
$\widetilde{x}(\theta)=(\sin(\theta),\cos(\theta)),$ и для
которой, в частности, выполнены условия леммы~\ref{MAlem} (см.
\cite{can}, формула~37). Подробное сопоставление с одной лишь
теоремой Мельникова в настоящей работе мотивировано еще и тем, что
недавно мы получили обобщение сформулированной выше теоремы
Малкина (см. \cite{nachr}, следствие~3.5), верное в пространстве
произвольной размерности. Подход, предложенный в \cite{nachr},
использует условие типа ($C_{MA}$) и к случаю Мельникова не
применим. Тоже самое условие предположено и в цитированной выше
работе \cite{anal}, где, в соответствующем случае, установлена
формула вида (\ref{Phi_mi}).

 В первом примере будет показано, что
теорема~\ref{t1_mi} может уточнять результат теоремы Мельникова, а
во втором, что она может устанавливать существование, по крайней
мере, двух близких циклу $\widetilde{x}$ $T$-периодических решений
в некоторых таких случаях, в которых теорема Мельникова
гарантирует существование, по крайней мере, одного.

Сравнение будет проводиться в случае главного резонанса, то есть,
 когда наименьший период цикла $\widetilde{x}$ совпадает с
наименьшим периодом возмущения, и более того, на его простом
подслучае, когда функция $M_E$ имеет ровно два нуля на интервале
$[0,T).$ Такое ограничение позволит применить теорему~\ref{analt}
и провести сравнение наиболее легко и наглядно.

\begin{exa}\hskip-0.2cm.\label{e4_mi} Пусть предложена система
\begin{equation}\label{mak_mi}
  \left(\begin{array}{c}
\dot x_1 \\ \dot
    x_2 \end{array}\right)=
    \left(\begin{array}{c}
     x_2(1-x_1^2-x_2^2)\\
     -x_1(1-x_1^2-x_2^2)
     \end{array}
    \right)+\varepsilon
    \left(\begin{array}{c}
     0 \\ \sin(wt) \end{array}\right).
\end{equation}
\end{exa}

\noindent При $\varepsilon=0$ система (\ref{mak_mi}) имеет вид
\begin{equation}\label{mak_np}
  \left(\begin{array}{c}
\dot x_1 \\ \dot
    x_2 \end{array}\right)=
    \left(\begin{array}{c}
     x_2(1-x_1^2-x_2^2)\\
     -x_1(1-x_1^2-x_2^2)
     \end{array}
    \right)
\end{equation}
и допускает семейство периодических орбит
$$
  \left\{\sqrt{1-\alpha}\left(\begin{array}{l} \sin(\alpha t) \\
  \cos(\alpha t)\end{array}\right)\right\}_{\alpha>0}.
$$
Рассмотрим задачу о существовании периодических орбит главного
резонанса, то есть задачу о возмущении орбиты
$$  \widetilde{x}(t)=\sqrt{1-w}\left(\begin{array}{l} \sin(wt) \\
  \cos(wt)\end{array}\right)
$$
периода $T=2\pi/w,$ совпадающего с периодом возмущения. Так как
$$
  \dot{\widetilde{x}}(t)=w\sqrt{1-w}\left(\begin{array}{c} \cos(wt)\\
  -\sin(wt)\end{array}\right),
$$
то функция Мельникова имеет вид
$$
  \widetilde{M}_E(\theta)=-\pi\sqrt{1-w}\sin(w\theta).
$$
Функция $\widetilde{M}_E$ имеет, очевидно, два простых нуля
$\theta_1=0$ и $\theta_2=\dfrac{\pi}{w}$ на интервале $[0,2\pi/w)$
и теорема Мельникова гарантирует, что {\it при любых $w\in(0,1)$ и
достаточно малых $\varepsilon>0$ система (\ref{mak_mi}) имеет по
крайней мере два $T$-периодических решения
$\widetilde{x}_{\varepsilon,1}$ и $\widetilde{x}_{\varepsilon,2},$
сходящихся при $\varepsilon\to 0$ к циклу $\widetilde{x}.$ }

Посмотрим теперь, какое утверждение позволяет получить
теорема~\ref{t1_mi}.

Заметим, что порождающая система (\ref{mak_np}) удовлетворяет
условиям леммы~\ref{symlem}, следовательно,
\begin{equation}\label{sledo}
M_E^s(\theta)= \widetilde{M}_E(\theta)\quad\mbox{и}\quad
  M_A^s(\theta)=-\int\limits_{s-T}^{s}{\rm
det}\left\|\left({\widehat{y}}(\tau),g(\tau-\theta,\widetilde{x}(\tau))\right)\right\|d\tau.
\end{equation}
Линеаризованная на $\widetilde{x}$ система (\ref{mak_np}) имеет
вид
\begin{equation}\label{mak_ls}
\left(\begin{array}{c} \dot y_1 \\ \dot y_2 \end{array}\right)=
\left(\begin{array}{cc} -2(1-w)\sin(wt)\cos(wt) & w-2(1-w)\cos^2(wt)\\
-w+2(1-w)\sin^2(wt) & 2(1-w)\sin(wt)\cos(wt)
\end{array}\right)\left(\begin{array}{c} y_1 \\ y_2
\end{array}\right).
\end{equation}
Решение $\widehat{y}$ системы (\ref{mak_ls}) с начальным условием
$$
  \widehat{y}(0)=\frac{1}{\|\dot{\widetilde{x}}(0)\|^2}\dot{\widetilde{x}}(0)^\bot=
  \frac{1}{w\sqrt{1-w}}
  \left(\begin{array}{c} 0 \\ 1 \end{array}\right)
$$
дается формулой
$$
  \widehat{y}(t)=\frac{1}{w\sqrt{1-w}}\left(\begin{array}{c} -2(1-w)t\cos(wt)+\sin(wt) \\
  2(1-w)t\sin(wt)+\cos(wt) \end{array}\right).
$$
Поэтому
\begin{eqnarray}
 M_A^s(0)&=&\frac{\pi}{w^3\sqrt{1-w}}\left(2(1-w)\sin^2(ws)-1\right),\nonumber\\
 M_A^s(\pi/w)&=&-\frac{\pi}{w^3\sqrt{1-w}}\left(2(1-w)\sin^2(ws)-1\right).\nonumber
\end{eqnarray}
Таким образом, для выполнения условия (\ref{L_cond})
теоремы~\ref{t1_mi} необходимо, чтобы
$$
  w\in\left(\frac{1}{2},1\right).
$$
При выполнении последнего условия имеем $M_A^s(0)M_A^s(\pi/w)<0$ и
в силу теоремы~\ref{analt}, условие (\ref{ind_mi})
теоремы~\ref{t1_mi} об индексе также выполнено. Итак, на основании
теоремы~\ref{t1_mi}, имеем: {\it при любых
$w\in\left(1/2,1\right)$  и достаточно малых $\varepsilon>0$
система (\ref{mak_mi}) имеет, по крайней мере, два
$2\pi/w$-периодических решения $\widetilde{x}_{\varepsilon,1}$ и
$\widetilde{x}_{\varepsilon,2},$ сходящихся при $\varepsilon\to 0$
к циклу $\widetilde{x}.$ Решение $\widetilde{x}_{\varepsilon,1}$
лежит строго внутри цикла $\widetilde{x},$ а решение
$\widetilde{x}_{\varepsilon,1}$ строго снаружи. Прочие
$2\pi/w$-периодические решения системы (\ref{mak_mi}) при
указанных $w$ и $\varepsilon$ также не пересекают цикл
$\widetilde{x}.$}

Мы получили, что для примера~\ref{e4_mi} область параметра $w>0,$
при которой применим метод Мельникова, шире области, при которой
применима теорема~\ref{t1_mi}, но в этой более узкой области
теорема~\ref{t1_mi} позволяет указать новые свойства периодических
решений главного резонанса.

Сейчас мы предъявим пример, в котором теорема~\ref{t1_mi}
указывает не только новые свойства периодических решений главного
резонанса, но и доказывает существование б\'ольшего их числа.

\begin{exa}\hskip-0.2cm.\label{e2_mi}
Действительно, подправим систему предыдущего примера следующим
образом
\begin{equation}\label{EX2}
  \begin{array}{lll}
    \dot x_1&=&x_2\left(1-\dfrac{1}{5}(x_1^2+x_2^2)\right),\\
    \dot
    x_2&=&-x_1\left(1-\dfrac{1}{5}(x_1^2+x_2^2)\right)+\varepsilon\left(\sin\left(\dfrac{4}{5}t\right)-x_1\right)^3+\varepsilon
    x_1
  \end{array}
\end{equation}
\end{exa}
и изучим возмущение порождающей орбиты
$$
  \widetilde{x}(t)=\left(\begin{array}{l} \sin\left(\dfrac{4}{5}t\right) \\
  \cos\left(\dfrac{4}{5}t\right)\end{array}\right)
$$
периода $T=\dfrac{5\pi}{2},$ совпадающего с периодом возмущения.
Соответствующая этой задаче функция Мельникова записывается как
$$
  M_E(\theta)=\frac{3\pi}{4}\sin\left(\frac{8}{5}t\right)-\frac{3\pi}{2}\sin\left(\frac{4}{5}t\right)
$$
и  допускает два нуля $\theta_1=0$ и $  \theta_2=\dfrac{5\pi}{4}.$
  Однако, только второй из них является простым, первый же
является кубическим, то есть $M'(0)=0,$ $M''(0)=0$ и
$M'''(0)=-\dfrac{288\pi}{125}.$

Таким образом, из теоремы Мельникова заключаем, что {\it при
достаточно малых $\varepsilon>0$ система (\ref{EX2}) имеет, по
крайней мере, одно $T$-периодическое решение
$\widetilde{x}_{\varepsilon},$ сходящееся при $\varepsilon\to 0$ к
циклу $\widetilde{x}.$ }

Легко проверить, что система (\ref{EX2}) с $\varepsilon=0,$
линеаризованная на $\widetilde{x},$ совпадает с системой
(\ref{mak_ls}), в которой взято $w=4/5.$ Поэтому, учитывая, что
 начальное условие решения $\widehat{y}$ дается
формулой
$$
  \widehat{y}(0)=\frac{1}{\|\dot{\widetilde{x}}(0)\|^2}\dot{\widetilde{x}}(0)^\bot=
  \frac{1}{4/5}
  \left(\begin{array}{c} 0 \\ 1 \end{array}\right),
$$
заключаем, что
$$
  \widehat{y}(t)=\frac{1}{4/5}\left(\begin{array}{c} -\dfrac{2}{5}t\cos\left(\dfrac{4}{5}t\right)+
  \sin\left(\dfrac{4}{5}t\right) \\
  \dfrac{2}{5}t\sin\left(\dfrac{4}{5}t\right)+\cos\left(\dfrac{4}{5}t\right) \end{array}\right).
$$
Пользуясь леммой~\ref{symlem}, а именно формулами (\ref{sledo}),
получаем для $M_A^s(0)$ и $M_A^s(5\pi/4)$ следующие выражения
\begin{eqnarray}
  M^s_A(0)&=&-\frac{25\pi}{64}\cos\left(\frac{8}{5}s\right)-\frac{25\pi}{16},\nonumber\\
  M^s_A(5\pi/4)&=&-\frac{25\pi}{64}\cos\left(\frac{16}{5}s\right)+\frac{75\pi}{64}\cos\left(\frac{8}{5}s\right)+
  \frac{125\pi}{16},\nonumber
\end{eqnarray}
для которых легко проверить, что $M^s_A(0)<0$ и $M^s_A(5\pi/4)>0$
при всех $s\in[0,5\pi/2].$

Значит, на основании теоремы~\ref{t1_mi} имеем: {\it при всех
достаточно малых $\varepsilon>0$ система (\ref{EX2}) имеет, по
крайней мере, два $T$-периодических решения
$\widetilde{x}_{\varepsilon,1}$ и $\widetilde{x}_{\varepsilon,2},$
сходящихся при $\varepsilon\to 0$ к циклу $\widetilde{x}.$ Решение
$\widetilde{x}_{\varepsilon,1}$ лежит строго внутри цикла
$\widetilde{x},$ а решение $\widetilde{x}_{\varepsilon,1}$ строго
снаружи. }

Таким образом, наличие у функции Мельникова кратных корней не
подрывает работоспособность теоремы~1.

Обратимся теперь к случаю, когда возмущение не является
дифференцируемой функцией. Отметим, что вопрос о существовании
периодических решений в возмущенных системах типа (\ref{ps_mi}) в
случае нулевой или линейной порождающей системы рассматривался
Ю.~А.~Митропольским \cite{mitro}, А.~М.~Самойленко \cite{sam},
Дж.~Мавеном \cite{maw1}, \cite{maw2}, А.~Буйка и Дж.~Либри
\cite{libri} и многими другими. В следующем примере
демонстрируется применение теоремы~\ref{t1_mi} к системам с
недифференцируемой правой частью при ненулевой нелинейной
порождающей системе. В качестве возмущения выбрана, так
называемая,  прыгающая нелинейность, см. \cite{lazer}. В качестве
порождающей системы выбрано уравнение Дуффинга.

\begin{exa}\hskip-0.2cm.\label{duff} Рассмотрим задачу о
существовании резонансных периодических решений у уравнения
Дуффинга с прыгающей нелинейностью
\begin{equation}\label{so}
  \ddot u+u+u^3=\varepsilon(\mu x^+_1+\nu
  x_1^-+\cos((1+\delta)t)).
\end{equation}
\end{exa}

Установим следующую простую лемму.

\begin{lem}\hskip-0.2cm.\label{jump}
Рассмотрим систему
\begin{equation}\label{lin}
  \begin{array}{lll}
    \dot x_1&=&x_2\\
    \dot x_2&=&-x_1+\varepsilon(\mu x^+_1+\nu x_1^-+\cos(t)),
  \end{array}
\end{equation}
где $a^+:=\max\{a,0\},$ $a^-:=\max\{-a,0\}.$
 Положим $\widetilde{\widetilde{x}}(t)=(\sin t,\cos t).$
 Тогда, если
$|\mu-\nu|\not=2,$ то соответствующие функции $M_A^s$ и $M_E^s$
удовлетворяют условию (\ref{L_cond}) теоремы~\ref{t1_mi}.
 Если же
$|\mu-\nu|<2,$ то $${\rm
ind}(\widetilde{\widetilde{x}},\Phi)\in\{0,2\}.$$
\end{lem}

\noindent {\rm Д о к а з а т е л ь с т в о.} Имеем
$$
  M^s_A(\theta)=\int\limits_0^{2\pi}\sin\tau\left(\mu \widetilde{x}^+_1(\tau)+\nu
  \widetilde{x}_1^-(\tau)+\cos(\tau-\theta)\right)d\tau=
$$
$$
  =\mu\int\limits_0^\pi\sin\tau\sin\tau
  d\tau-\nu\int\limits_\pi^{2\pi}\sin\tau\sin\tau
  d\tau+\int\limits_0^{2\pi}\sin\tau\cos\tau
  d\tau\cos\theta+\int\limits_0^{2\pi}\sin\tau\sin\tau d\tau\sin\theta=
$$
$$
  = \mu\frac{\pi}{2}-\nu\frac{\pi}{2}+\pi\sin\theta,
$$
$$
  M_E^s(\theta)=\int\limits_0^{2\pi}\cos\tau\left(\mu \widetilde{x}^+_1(\tau)+\nu
  \widetilde{x}_1^-(\tau)+\cos(\tau-\theta)\right)d\tau=
$$
$$
  =\mu\int\limits_0^\pi\cos\tau\sin\tau
  d\tau-\nu\int\limits_\pi^{2\pi}\cos\tau\sin\tau
  d\tau+\int\limits_0^{2\pi}\cos\tau\cos\tau
  d\tau\cos\theta+\int\limits_0^{2\pi}\cos\tau\sin\tau d\tau\sin\theta=
$$
$$
  =\pi\cos\theta.
$$
Так как по условию леммы $|\mu-\nu|<2,$ то
$\theta_0=\arcsin\dfrac{-\mu+\nu}{2}$ будет единственным нулем
функции $M_E$ на интервале $[-\pi/2,\pi/2].$ Поэтому на интервале
$[0,2\pi)$ функция $M_E$ имеет ровно два корня
$$
\begin{array}{l}
  \theta_1=\left\{\begin{array}{l}\theta_1=\theta_0,\quad\mbox{если
  }\theta_0\ge 0,\\
  \theta_1=\theta_0+\pi,\quad\mbox{в противном
  случае,}\end{array}\right.\\
  \theta_2=\theta_1+\pi.
\end{array}
$$
Функция $M_E$ имеет ровно два нуля $\theta_1=\pi/2$ и
$\theta_2=3\pi/2$ на интервале $[0,2\pi).$ Для функции
$M_A(\theta)$ в этих точках имеем
$$
  M_A(\pi/2)=\mu\frac{\pi}{2}-\nu\frac{\pi}{2}+\pi,\qquad
  M_A(3\pi/2)=\mu\frac{\pi}{2}-\nu\frac{\pi}{2}-\pi.
$$
Таким образом, если $|\mu-\nu|\not=2,$ то
$$
  M_A(\pi/2)\not=0,\qquad
  M_A(3\pi/2)\not=0.
$$
Если же $|\mu-\nu|<2,$ то
$$
  M_A(\pi/2)>0,\qquad
  M_A(3\pi/2)<0
$$
и в силу теоремы~\ref{analt} имеем ${\rm
ind}(\widetilde{\widetilde{x}},\Phi)\in\{0,2\}.$

Лемма доказана.

Итак, вернемся к системе (\ref{so}). Обозначим через $u_\delta$
единственную с точностью до сдвига периодическую орбиту
порождающего уравнения
$$
  \ddot u+u+u^3=0
$$
  с наименьшим периодом $\dfrac{2\pi}{1+\delta}.$ Если $u$ -- решение уравнения (\ref{so}),
  то
$v=(u,\dot u)$  удовлетворяет системе
\begin{equation}\label{duf}
  \begin{array}{lll}
    \dot v_1&=&v_2\\
    \dot v_2&=&-v_1- v_1^3+\varepsilon(\mu v^+_1+\nu v_1^-+\cos((1+\delta)t)),
  \end{array}
\end{equation}
обратно, если $v$ -- решение системы (\ref{duf}), то $v_1$ --
решение системы (\ref{so}). Без ограничения общности для
предлагаемого ниже утверждения можно считать, что $\dot
u_\delta(0)=0$ и $u_\delta(0)>0.$ Заменой переменных
$$
  x(t)=\frac{v(t)}{u_\delta(0)}
$$
перейдем от системы (\ref{duf}) к системе
\begin{equation}\label{duf1}
  \begin{array}{lll}
    \dot x_1&=&x_2\\
    \dot x_2&=&-x_1-(u_\delta(0))^2x_1^3+\varepsilon
    \dfrac{1}{u_\delta(0)}\cos((1+\delta)t) +\varepsilon(\mu x^+_1+\nu x_1^-).
  \end{array}
\end{equation}
Пусть $|\mu-\nu|<2.$ Установим, что существует $\delta_0>0$ такое,
что при $\delta\in(0,\delta_0]$ условия теоремы~\ref{t1_mi},
связанные с функциями $M_E^s,$ $M_A^s$ и $\Phi,$ для системы
(\ref{duf1}) выполнены с
$\widetilde{x}(t)=\dfrac{v_\delta(t)}{u_\delta(0)}$ и
$T=\dfrac{2\pi}{1+\delta}.$ Для этого, в свою очередь, достаточно
установить аналогичное утверждение для системы
\begin{equation}\label{duf11}
  \begin{array}{lll}
    \dot x_1&=&x_2\\
    \dot x_2&=&-x_1-(u_\delta(0))^2x_1^3+\varepsilon
    \cos((1+\delta)t)+\varepsilon(\mu x^+_1+\nu x_1^-).
  \end{array}
\end{equation}
Период орбит порождающей системы
\begin{equation}\nonumber
  \begin{array}{lll}
    \dot x_1&=&x_2\\
    \dot x_2&=&-x_1-(u_\delta(0))^2x_1^3
  \end{array}
\end{equation}
изменяется монотонно от $2\pi$ до $0,$ когда начальное условие
орбиты изменяется от $(0,0)$ до $(+\infty,0)$ (см., напр.,
\cite[пример с.~250]{guck}). Следовательно, $u_\delta(0)\to 0,$
когда $\delta\to 0.$ Но для $\delta=0$ справедливость желаемого
для системы (\ref{duf11}) утверждения следует из леммы~\ref{jump},
следовательно, это утверждение остается справедливым и при малых
$\delta>0.$

Итак, установлено, что {\it если $|\mu-\nu|<2,$ то существует
$\delta_0>0$ такое, что каждому $\delta\in[0,\delta_0]$
соответствует $\varepsilon_0>0$ такое, что:

1) при каждом $\varepsilon\in(0,\varepsilon_0]$ возмущенное
уравнение Дуффинга (\ref{so}) имеет, по крайней мере, два
$\frac{2\pi}{1+\delta}$-периодических решения
$\widetilde{u}_{\delta,\varepsilon},$
$\widetilde{\widetilde{u}}_{\delta,\varepsilon}$ таких, что кривая
$t\to\left(\widetilde{u}_{\delta,\varepsilon}(t),\dot{\widetilde{u}}_{\delta,\varepsilon}(t)\right)$
лежит строго внутри кривой
$t\to\left(u_\delta(t),\dot{{u}}_\delta(t)\right),$ а кривая
$t\to\left(\widetilde{\widetilde{u}}_{\delta,\varepsilon}(t),\dot{\widetilde{\widetilde{u}}}_{\delta,\varepsilon}
(t)\right)$ строго снаружи;

2) для всякого $\frac{2\pi}{1+\delta}$-периодического решения $u$
системы (\ref{so}) с $\varepsilon\in(0,\varepsilon_0]$ кривая
$t\to(u(t),\dot u(t))$ не имеет точек пересечения с  кривой
$t\to\left(u_\delta(t),\dot{{u}}_\delta(t)\right);$

3) решения $\widetilde{u}_{\delta,\varepsilon}$  и
$\widetilde{\widetilde{u}}_{\delta,\varepsilon}$ удовлетворяют
условию
$$
  \widetilde{u}_{\delta,\varepsilon}(t)\to u_\delta\left(t-\widetilde{\theta}\right)
  \quad\mbox{и}\quad
  \widetilde{\widetilde{u}}_{\delta,\varepsilon}(t)\to u_\delta\left(t-\widetilde{\widetilde{\theta}}\right)
  \quad\mbox{при}\ \varepsilon\to 0
$$
для некоторых
$\widetilde{\theta},\widetilde{\widetilde{\theta}}\in\left[0,\frac{2\pi}{1+\delta}\right].$
}

Если $\nu=\mu=0,$ то полученное утверждение является добавлением к
известным результатам о качественном поведении периодических
траекторий уравнения Дуффинга, см., напр., А.~Д.~Морозов
\cite{mord}.

\subsection{Вырожденные резонансы. Сопоставление с теоремой
Йагасаки}\label{degres}
\subsubsection{}\vskip-1cm
\addtocounter{sect}{1}

В настоящем параграфе рассматривается случай, когда
$T$-периодический цикл $\widetilde{x}$ является вырожденным, то
есть все решения линеаризованной системы (\ref{ls_mi}) являются
$T$-периодическими. $T$-периодические решения возмущенной системы,
порожденные такими циклами, называются вырожденными резонансами.
Если $\widetilde{x}$ вложен в семейство циклов
$\left\{\widetilde{x}_\alpha\right\}_{\alpha>0}$ автономной
системы с периодами $T(\alpha),$ то есть, если
$$
  \widetilde{x}=\widetilde{x}_{\alpha_0}
$$
для некоторого $\alpha_0>0,$ то, как показывает нижеследующая
лемма, предположение о вырожденности почти всегда выполнено в
случае, когда $T(\alpha_0)$ является критическим периодом, то есть
\begin{equation}\label{deg}
  T'(\alpha_0)=0.
\end{equation}

Отметим, что задача о существовании критических периодов для
циклов, вложенных в семейство циклов автономной системы,
интенсивно исследуется, см. \cite{crit1}, \cite{crit2},
\cite{crit3}, \cite{yulin}.

Не ограничивая общности можем считать далее, что
$$
  \widetilde{x}_\alpha(0)=\left(\begin{array}{c}0\\
  J(\alpha)\end{array}\right).
$$

\begin{lem}\hskip-0.2cm.\label{krit} Если для орбиты $\widetilde{x}_{\alpha_0}$ выполнены
условия  $T'(\alpha_0)=0,$ $J'(\alpha_0)\not=0$ и первая
компонента вектора $\dot{\widetilde{x}}_{\alpha_0}(0)$ отлична от
нуля, то каждое решение линеаризованной системы (\ref{ls_mi}) с
$\widetilde{x}=\widetilde{x}_{\alpha_0}$ является
$T(\alpha_0)$-периодическим.
\end{lem}

\noindent {\rm Д о к а з а т е л ь с т в о.} Нам будет удобно
использовать следующее обозначение
$$
  x(t,\alpha):=\widetilde{x}_\alpha(t).
$$
 Так как правая часть порождающей системы (\ref{np_mi})
непрерывно дифференцируема, то (см., напр., Понтрягин \cite{pont},
Гл.~4, \S~24, теорема~17]) функция $(t,\alpha)\to x(t,\alpha)$
непрерывно дифференцируема по совокупности переменных.
Дифференцируя тождество
$$
x'_t(t,\alpha)=f(x(t,\alpha))
$$
по $\alpha,$ получаем
$$x'_t {}'_\alpha(t,\alpha)=f'(x(t,\alpha))x'_\alpha,$$
следовательно, $\widetilde{y}=x'_\alpha(\cdot,\alpha_0)$ является
решением линеаризованной системы (\ref{ls_mi}), в которой
$\widetilde{x}=x(\cdot,\alpha_0).$ Имеем
\begin{eqnarray}
  \widetilde{y}(T(\alpha_0))-\widetilde{y}(0)&=&\lim_{\Delta\to 0}\frac{x(T(\alpha_0),\alpha_0+\Delta)-
  x(T(\alpha_0),\alpha_0)}{\Delta}-\lim_{\Delta\to
  0}\frac{x(0,\alpha_0+\Delta)-x(0,\alpha_0)}{\Delta}=\nonumber\\
  & =&\lim_{\Delta\to 0}\frac{x(T(\alpha_0),\alpha_0+\Delta)-x(0,\alpha_0+\Delta)-
  x(T(\alpha_0),\alpha_0)+x(0,\alpha_0)}{\Delta}=\nonumber\\
  &=&\lim_{\Delta\to
  0}\frac{x(T(\alpha_0),\alpha_0+\Delta)-x(0,\alpha_0+\Delta)}{\Delta}=\nonumber\\
  &=&\left(x(T(\cdot),\alpha_0)\right)'(\alpha_0)=x'_t(T(\alpha_0),\alpha_0)T'(\alpha_0)=0,
\end{eqnarray}
то есть $\widetilde{y}$ является $T(\alpha_0)$-периодическим
решением системы (\ref{ls_mi}). Но $x(0,\alpha)=\left(\begin{array}{c}0\\
  J(\alpha)\end{array}\right),$ следовательно
$$
  \widetilde{y}(0)=x'_\alpha(0,\alpha_0)=\left(\begin{array}{c}0\\
  J'(\alpha_0)\end{array}\right).
$$
В тоже время, согласно условиям леммы, первая компонента вектора
$\dot{\widetilde{x}}_{\alpha_0}(0)$ отлична от нуля, значит
 $\widetilde{y}$ и $x'_t(\cdot,\alpha_0)$ -- два
линейно-независимых $T(\alpha_0)$-периодических решения
(двумерной) системы (\ref{ls_mi}), что, очевидно, влечет желаемое
утверждение.

Лемма доказана.

\vskip0.3cm

Лемма~\ref{krit} является критерием вырожденности цикла,
вложенного в семейство циклов. Однако, для вырожденности цикла
вовсе не необходимо чтобы он был вложен в семейство циклов,
вырожденными могут быть и изолированные циклы. Для изучения
существования в возмущенной системе (\ref{ps_mi}) вырожденных
резонансов в этом последнем случае может, вообще говоря,
использоваться общая теорема Рума-Чиконе (\cite{chic},
теорема~4.1), но она работает только в случае, когда возмущение
зависит от фазовой переменной (см. \cite{chic}, формула~2.7), что
не требуется в предлагаемых ниже теоремах.

Если цикл $\widetilde{x}$ является вырожденным, то функции $M_A^s$
и $M_E^s$ очевидно, не зависят от $s$ и теорема~\ref{t1_mi}
принимает следующий вид.

\begin{thm}\hskip-0.2cm.\label{t1d_mi}
 Пусть для вырожденного цикла $\widetilde{x}$ выполнено  условие ($C$) и возмущение
непрерывно. Предположим, что для любого $\theta_0\in[0,T]$ такого,
что $M_E(\theta_0)=0$ имеем
$$
  M_A(\theta_0)\not=0.
$$
Тогда при достаточно малых $\varepsilon>0$ всякое
$T$-периодическое решение $\widetilde{x}_\varepsilon$ возмущенной
системы (\ref{ps_mi}) необходимо таково, что
\begin{equation}\label{us2d_mi}
  \widetilde{x}_\varepsilon(t)\not=\widetilde{x}(s)\quad\mbox{при всех
  }t,s\in[0,T].
\end{equation}
Если же имеем еще и
$$
{\rm ind}(\widetilde{x},\Phi)\not=1,
$$
то при всех достаточно малых $\varepsilon>0$ возмущенная система
(\ref{ps_mi}) действительно имеет, по крайней мере, два
$T$-периодических решения $\widetilde{x}_{\varepsilon,1}$ и
$\widetilde{x}_{\varepsilon,2},$ удовлетворяющих (\ref{us2d_mi}).
Оба решения сходятся к $\widetilde{x}$ при $\varepsilon\to 0.$
Кроме того, решение $\widetilde{x}_{\varepsilon,1}$ лежит  внутри
цикла $\widetilde{x},$ а решение $\widetilde{x}_{\varepsilon,2}$
снаружи.
\end{thm}

Условия теоремы~\ref{t1_mi} приобретают максимально простой вид,
если дополнительно к вырожденности известно, что функция $M_E$
имеет ровно два нуля на $[0,T).$ По этой причине мы сформулируем
соответствующее утверждение как отдельную теорему.

\begin{thm}\hskip-0.2cm.\label{t2_mi}
 Пусть для вырожденного цикла $\widetilde{x}$ выполнено  условие ($C$) и возмущение
непрерывно. Предположим, что функция $M_E$ имеет ровно два нуля
$\theta_1$ и $\theta_2$ на интервале $[0,T)$ и
$$
  M_A(\theta_1)\cdot M_A(\theta_2)\not=0.
$$
Тогда при достаточно малых $\varepsilon>0$ всякое
$T$-периодическое решение $\widetilde{x}_\varepsilon$ возмущенной
системы (\ref{ps_mi}) необходимо удовлетворяет условию
(\ref{us2d_mi}). Если же дополнительно известно, что
$$
  M_A(\theta_1)\cdot M_A(\theta_2)<0,
$$
то при всех достаточно малых $\varepsilon>0$ возмущенная система
(\ref{ps_mi}) действительно имеет, по крайней мере, два
$T$-периодических решения $\widetilde{x}_{\varepsilon,1}$ и
$\widetilde{x}_{\varepsilon,2},$ удовлетворяющих (\ref{us2d_mi}).
Оба решения сходятся к $\widetilde{x}$ при $\varepsilon\to 0.$
Кроме того, решение $\widetilde{x}_{\varepsilon,1}$ лежит  внутри
цикла $\widetilde{x},$ а решение $\widetilde{x}_{\varepsilon,2}$
снаружи.
\end{thm}

Ясно, что верна теорема, полученная из теоремы~\ref{t2_mi}
перестановкой функций $M_E$ и $M_A$ между собой.

Насколько известно автору, в литературе исследован только тот
случай, когда вырожденный цикл вложен в семейство циклов и имеет
критический период. Поэтому, именно такую ситуацию мы предположим
далее для тестирования теоремы~\ref{t2_mi}. В этом случае условия,
связанные с применением теорем о неявной функции, напротив,
усложняются и существенно зависят от того, какова кратность нуля
$\alpha_0$ для функции $T.$ Так К.~Йагасаки \cite{yaga} рассмотрел
случай, когда
$$
  T'(\alpha_0)=0, \ \ T''(\alpha_0)\not=0.
$$
Помимо условия существования простого корня $\theta_0$ у функции
Мельникова, в работе цитированного автора требуется, чтобы
некоторая вспомогательная функция $N,$ формула которой содержит
двойной интеграл, имела в точке $\theta_0$ определенный знак.
Полученную теорему Йагасаки иллюстрирует на нескольких примерах,
некоторую модификацию одного из них (\cite{yaga}, пример~\S~6.4)
мы выберем для сравнения этой теоремы с теоремой~\ref{t2_mi}.

\begin{exa}\hskip-0.2cm. Действительно, рассмотрим систему
\begin{equation}\label{yag}
  \begin{array}{lll}
    \dot x_1&=&x_2\left(\dfrac{1}{4}(x_1^2+x_2^2-2)^p+1\right)\\
    \dot
    x_2&=&-x_1\left(\dfrac{1}{4}(x_1^2+x_2^2-2)^p+1\right)+\varepsilon\sin(t).
  \end{array}
\end{equation}
\end{exa}
Порождающая система допускает семейство циклов
$$
  \widetilde{x}_\alpha(t)=\left(\begin{array}{c}
  \alpha\sin\left(\dfrac{2\pi}{T(\alpha)}t\right)\\
  \alpha\cos\left(\dfrac{2\pi}{T(\alpha)}t\right)
  \end{array}\right)
$$
 с периодами
 $$
   T(\alpha)=\frac{2\pi}{\dfrac{1}{4}(\alpha^2-2)^p+1}.
 $$
 Исследуем возмущение цикла, соответствующего $\alpha=\sqrt{2},$
 то есть цикла
 $$
   \widetilde{x}(t)=\left(\begin{array}{c}
  \sqrt{2}\sin t\\
  \sqrt{2}\cos t
  \end{array}\right),
 $$
 для которого
 $$
   T\left(\sqrt{2}\right)=2\pi,\ \ T'\left(\sqrt{2}\right)=...=T^{(p-1)}\left(\sqrt{2}\right)=0,
   \ \ T^{(p)}\left(\sqrt{2}\right)\not=0.
 $$

В случае $p=2$ для системы (\ref{yag}) Йагасаки вычисляет
упомянутую функцию $N$ и приводит следующее утверждение (см.
\cite{yaga}, теорема~6.4):  {\it Для всех достаточно малых
$\varepsilon>0$ система (\ref{yag}) имеет, по крайней мере, два
$T$-периодических решения $\widetilde{x}_{\varepsilon,1}$ и
$\widetilde{x}_{\varepsilon,2},$ сходящихся к $\widetilde{x}$ при
$\varepsilon\to 0.$ }

Так как при $p\not=2$ критический период $2\pi$ не является
двукратным, то при $p\not=2$ результат Йагасаки не применим.

Попробуем использовать теорему~\ref{t2_mi}. Линеаризованная на
цикле $\widetilde{x}$ порождающая система при любом натуральном
$p$ имеет вид
$$
  \left(\begin{array}{c} \dot y_1\\
  \dot y_2\end{array}\right)=\left(\begin{array}{cc} 0 & 1\\ -1 &
  0\end{array}\right)\left(\begin{array}{c} y_1\\
  y_2\end{array}\right)
$$
и, следовательно,
$$
   \widehat{y}(t)=\frac{1}{\sqrt{2}}\left(\begin{array}{c}
  \sin t\\
  \cos t
  \end{array}\right).
$$
 Таким образом, при любом натуральном $p$ имеем
$$
  M_E(\theta)=-\sqrt{2}\pi\sin(\theta),\ \ \
  M_A(\theta)=-\frac{1}{\sqrt{2}}\pi\cos(\theta).
$$
 Функция $M$ имеет два нуля $\theta_1=0$ и
$\theta_2=\pi$ на интервале $[0,T),$ причем $M_A(0)\cdot
M_A(\pi)<0.$ Таким образом, выполнены условия теоремы~\ref{t2_mi},
на основании которой получаем утверждение:
 {\it Пусть $p\in\mathbb{N}$ -- произвольное число.
 Для всех достаточно малых $\varepsilon>0$
система (\ref{yag}) имеет, по крайней мере, два $T$-периодических
решения $\widetilde{x}_{\varepsilon,1}$ и
$\widetilde{x}_{\varepsilon,2},$ сходящихся к $\widetilde{x}$ при
$\varepsilon\to 0.$ Решения $\widetilde{x}_{\varepsilon,1}$ лежат
строго внутри цикла $\widetilde{x}$ а решения
$\widetilde{x}_{\varepsilon,1}$ строго снаружи. }

Таким образом, при $p=2$ полученное из теоремы~\ref{t2_mi}
утверждение уточняет результат  Йагасаки и дает точно такое
заключение при любой другой кратности вырождения, где указанный
результат не применим. В тоже время стоит отметить, что
рассматриваемая теорема Йагасаки (\cite{yaga}, теорема~6.4)
гарантирует для некоторых систем (при $p=2$)   существование, по
крайней мере, четырех периодических решений, в то время как
теорема~\ref{t2_mi} всегда гарантирует существование, по крайней
мере, двух.

В заключение параграфа отметим, что поведение траекторий
возмущенной гамильтоновой системы вблизи цикла с критическим
периодом исследовано в \cite{mor} и \cite{shil}.

\subsection{Расположение устойчивых и неустойчивых периодических решений}
\subsubsection{}\vskip-1cm
\addtocounter{sect}{1}

Теоремы Малкина, Мельникова и Йагасаки, как основанные на теореме
о неявной функции, предоставляют также информацию об устойчивости
рожденных из цикла $\widetilde{x}$ периодических решений
возмущенной системы (\ref{ps_mi}). В настоящем параграфе будет
показано, что знание индекса  ${\rm ind}(\widetilde{x},\Phi)$
позволяет в некоторых случаях утверждать, по какую сторону от
цикла находится, по крайней мере, одно устойчивое решение, а по
какую сторону, по крайней мере, одно неустойчивое.

Предположим, что

$(A_{\mathcal{P}})$ решение $x=:\Omega_\varepsilon(\cdot,t_0,\xi)$
возмущенной системы (\ref{ps_mi}) с начальным условием
$x(t_0)=\xi$ существует, единственно и продолжимо на отрезок
$[0,T]$ при любых $t_0\in[0,T],$ $\xi\in\mathbb{R}^2$ и
$\varepsilon>0.$

Cледовательно, для системы (\ref{ps_mi}) при любых $\varepsilon>0$
определен оператор Пуанкаре-Андронова
$\mathcal{P}_\varepsilon=\Omega_\varepsilon(T,0,\cdot),$
соответствующий задаче о $T$-периодических решениях для
(\ref{ps_mi}).

\begin{thm}\hskip-0.2cm.\label{andr}
 Пусть возмущение в (\ref{ps_mi})
непрерывно. Пусть выполнено условие (\ref{L_cond})
теоремы~\ref{t1_mi} и условие $(A_{\mathcal{P}}).$
 \noindent Тогда существует $\varepsilon_0>0$
такое, что
$$
  \widetilde{x}(\theta)\not=\mathcal{P}_\varepsilon(\widetilde{x}(\theta))\quad\mbox{при
  всех }\  \theta\in[0,T],\ \varepsilon\in(0,\varepsilon_0]
$$
и
$$
  {\rm ind}(\widetilde{x},I-\mathcal{P}_\varepsilon)={\rm ind}(\widetilde{x},\Phi),\qquad\varepsilon\in(0,\varepsilon_0].
$$
\end{thm}

\noindent{\rm Д о к а з а т е л ь с т в о.} Пусть $Q_\varepsilon$
-- интегральный оператор, соответствующий задаче о
$T$-периодических решениях для возмущенной системы (\ref{ps_mi}) и
определяемый формулой (\ref{operator}). Обозначим через
$U\subset\mathbb{R}^2$ внутренность цикла $\widetilde{x}$ и
положим
$$
  W_\varepsilon=\left\{\hat x:
  C([0,T],\mathbb{R}^2):\Omega_\varepsilon(0,t,\hat x(t))\in
  U,\ {\rm for\ any\ }t\in[0,T]\right\}.
$$
 Покажем, что существует
$\varepsilon_0>0$ такое, что для любого
$\varepsilon\in(0,\varepsilon_0]$ и любого
$\alpha\in[0,\varepsilon_0]$ выполнено:
\begin{equation}\label{ggg}
  \mbox{\rm если}\ Q_\varepsilon x=x\ \mbox{\rm и}\ x\in \overline{W}_\alpha\ \mbox{\rm то,}\
  x\in W_0.
\end{equation}

Предположим противное, следовательно существуют последовательности
$\{\varepsilon_n\}_{n\in\mathbb{N}}\subset(0,\varepsilon_0],$
$\varepsilon_n\to 0$ при $n\to\infty,$
$\{\alpha_n\}_{n\in\mathbb{N}}\subset(0,\varepsilon_0],$
$\{x_n\}_{n\in\mathbb{N}}\subset C([0,T],\mathbb{R}^2),$ $x_n\in
\overline{W}_{\alpha_n}$ такие, что $Q_{\varepsilon_n}x_n=x_n$ и
$x_n\not\in W_0.$ Так как $x_n\in \overline{W}_{\varepsilon_n},$
то $x_n(0)\in U.$ С другой стороны из соотношения $x_n\not\in W_0$
заключаем, что при любом $n\in\mathbb{N}$ существует $t_n\in(0,T]$
такое, что $x_n(t_n)\in\partial U,$ в чем противоречие с
утверждением~1) теоремы~\ref{fromnon}.

 Из (\ref{ggg}) и утверждения~1) теоремы~\ref{fromnon} заключаем, что степень $d(I-Q_\varepsilon,W_\alpha)$ определена при
 любом  $\alpha\in
(0,\varepsilon_0]$ и
$$
  d(I-Q_\varepsilon,W_\varepsilon)=d(I-Q_\varepsilon,W_0),\quad\varepsilon\in(0,\varepsilon_0].
$$
Из принципа родственности (см. \cite{krazab}, теорема~28.5)
следует, что
$$
  d(I-Q_\varepsilon,W_\varepsilon)={\rm ind}(\widetilde{x},I-\mathcal{P}_\varepsilon).
$$
С другой стороны, в силу утверждения~2) теоремы~\ref{fromnon}
имеем
$$
  d(I-Q_\varepsilon,W_0)={\rm ind}(\widetilde{x},\Phi).
$$

Теорема доказана.

Нам также необходима следующая лемма, установленная Капетто,
Мавеном и Занолином.

{\rm Л е м м а \ \ К а п е т т о - М а в е н а - З а н о л и н а.
}\ \  (см. \cite{camaza}, Следствие~2). {\it Пусть выполнено
условие ($C$)  и возмущение непрерывно. Пусть выполнено условие
$(A_{\mathcal{P}}).$ Тогда для любой достаточно близкой к циклу
$\widetilde{x}$ и непересекающей его $T$-периодической функции
$\widehat{x}$ имеет место равенство
\begin{equation}\label{rrr}
  {\rm ind}(\widehat{x},I-\mathcal{P}_0)=1.
\end{equation}}
Доказательство леммы Капетто-Мавена-Занолина использует  теорему
Купки-Смейла (см. \cite{palis}, Гл.~3).

Итак, пользуясь индексом ${\rm ind}(\widetilde{x},\Phi),$ имеем
следующую информацию о  расположении устойчивых и неустойчивых
$T$-периодических решений возмущенной системы (\ref{ps_mi}) вблизи
порождающего цикла $\widetilde{x}.$

\begin{thm}\label{stab}\hskip-0.2cm. Пусть выполнены условия теоремы (\ref{andr}), а также
условие ($C$).  Выберем произвольные $T$-периодические функции
$\widetilde{x}^-$ и $\widetilde{x}^+$ такие, что $\widetilde{x}^-$
лежит внутри цикла $\widetilde{x},$ а $\widetilde{x}^+$ --
снаружи, причем в области $V^-,$ заключенной между
$\widetilde{x}^-$ и $\widetilde{x},$ а также в области $V^+,$
заключенной между $\widetilde{x}$ и $\widetilde{x}^+,$ нет
$T$-периодических решений порождающей системы (\ref{np_mi}). Пусть
дополнительно известно, что все неподвижные точки оператора
Пуанкаре-Андронова $\mathcal{P}_\varepsilon$ в $V^-$ и $V^+$
являются простыми. Тогда, если
$$
  {\rm ind}(\widetilde{x},\Phi)>1\quad \left({\rm
  ind}(\widetilde{x},\Phi)>1\right),
$$
то множество $V^+$ $\left(V^-\right)$ содержит, по крайней мере,
$$\mu=\left|{\rm ind}(\widetilde{x},\Phi)-1\right|$$
неподвижных точек оператора $\mathcal{P}_\varepsilon,$ являющихся
седлами, а множество $V^-$ $\left(V^+\right)$ содержит, по крайней
мере, $\mu$ неподвижных точек оператора $\mathcal{P}_\varepsilon,$
каждая из которых либо узел, либо фокус.
\end{thm}

Отметим, что выбор указанных в формулировке циклов
$\widetilde{x}^-$ и $\widetilde{x}^+$ возможен в силу
предположения ($C$).

\noindent {\rm Д о к а з а т е л ь с т в о. } Пусть
$$
{\rm ind}(\widetilde{x},\Phi)>1.
$$
Тогда, в силу теоремы~\ref{andr} и теоремы
Капетто-Мавена-Занолина, имеем
$$
  d(\Phi,V^-)>\mu\quad\mbox{и}\quad d(\Phi,V^+)
  <\mu.
$$
Следовательно (см. \cite{andr}, Гл.~V, \S~11, теорема~26 и
лемма~1), существует, по крайней мере, $\mu$ точек
$\xi^-_1,...,\xi^-_\mu,$
$$
\xi^-_i\in V^-\ \mbox{для любого } i\in\overline{1,\mu},
$$
и, по крайней мере, $\mu$ точек $\xi^+_1,...,\xi^+_\mu,$
$$
\xi^+_i\in V^+\ \mbox{для любого } i\in\overline{1,\mu},
$$
 таких, что
$$
  {\rm ind}({\xi}^-_i)=+1\quad\mbox{и}\quad {\rm
  ind}({\xi}^+_i)=-1\ \mbox{для любого } i\in\overline{1,\mu},
$$
то есть (см. \cite{andr}, Гл.~V, \S~11, теорема~30) каждая из
точек ${\xi}^-_i,$  $i\in\overline{1,\mu},$ является узлом или
фокусом, а каждая из точек $\xi^+_i,$  $i\in\overline{1,\mu},$ --
седлом.

Случай, когда $ {\rm ind}(\widetilde{x},\Phi)<1,$ рассматривается
аналогично.

Теорема доказана.

Развитие теоремы~\ref{stab} может быть получено на основе схем,
предложенных в \cite{ortega}.

\subsection{Заключение}\nonumber
\subsubsection{}\vskip-1cm
\addtocounter{sect}{1}

Итак, в работе для исследования задачи о рождении
$T$-периодических решений возмущенной системы (\ref{ps_mi}) из
$T$-периодического цикла $\widetilde{x}$ порождающей системы
(\ref{np_mi}) предложена новая характеристика порождающего цикла
${\rm ind}(\widetilde{x},\Phi).$ Даны условия, при которых из
свойства ${\rm ind}(\widetilde{x},\Phi)\not=1$ следует, что цикл
$\widetilde{x}$ порождает по крайней мере два $T$-периодических
решения возмущенной системы (\ref{ps_mi}), лежащих по разные
стороны от $\widetilde{x}.$ В зависимости от того ${\rm
ind}(\widetilde{x},\Phi)>1$ или ${\rm ind}(\widetilde{x},\Phi)<1$
делаются некоторые выводы о том, по какую сторону от
$\widetilde{x}$ рождаются устойчивые $T$-периодические решения и
по какую неустойчивые.

\addcontentsline{toc}{abcd}{Литература}

\begin{center}
\bf{ Литература \vskip-2cm $\ $}
\end{center}

\def\refname{}

\end{document}